\documentclass{article}

\usepackage{arxiv}

\usepackage[utf8]{inputenc} 
\usepackage[T1]{fontenc}    
\usepackage{hyperref}       
\usepackage{url}            
\usepackage{booktabs}       
\usepackage{amsfonts}       
\usepackage{nicefrac}       
\usepackage{microtype}      
\usepackage{lipsum}
\usepackage{graphicx}

\usepackage{natbib}
\usepackage{algorithmic}
\usepackage{algorithm}
\usepackage{amsmath,amssymb}
\usepackage{amsthm}
\usepackage{bm}
\usepackage{mathrsfs}
\usepackage{subcaption}

\newcommand{\argmin}{\operatornamewithlimits{argmin}}

\title{Multi-stage stochastic linear programming for shared autonomous vehicle system operation and design with on-demand and pre-booked requests}

\author{
 Riki Kawase \thanks{Corresponding author. Submitted to TRB Annual Meeting 2025} \\
  Tokyo Institute of Technology \\
  \texttt{kawase.r.ac@m.titech.ac.jp} \\
}

\begin{document}
\maketitle
\begin{abstract}
This study presents optimization problems to jointly determine long-term network design, mid-term fleet sizing strategy, and short-term routing and ridesharing matching in shared autonomous vehicle (SAV) systems with pre-booked and on-demand trip requests.
Based on the dynamic traffic assignment framework, multi-stage stochastic linear programming is formulated for joint optimization of SAV system design and operations.
Leveraging the linearity of the proposed problem, we can tackle the computational complexity due to multiple objectives and dynamic stochasticity through the weighted sum method and stochastic dual dynamic programming (SDDP).
Our numerical examples verify that the solution to the proposed problem obtained through SDDP is close enough to the optimal solution.
We also demonstrate the effect of introducing pre-booking options on optimized infrastructure planning and fleet sizing strategies.
Furthermore, dedicated vehicles to pick-up and drop-off only pre-booked travelers can lead to incentives to reserve in advance instead of on-demand requests with little reduction in system performance.
\end{abstract}


\section{Introduction}
Shared autonomous vehicle (SAV) systems will emerge as an innovative solution to urban traffic congestion in the future \citep{narayanan2020shared}.
With the recent advances in autonomous driving technologies and the ubiquity of the Internet, autonomous vehicles shared by society can provide optimized routes and ridesharing matching for on-demand requests from users via mobile applications \citep{levin2017general}.
Historical request data can help operators acquire statistical information on spatio-temporal traffic demand, thereby addressing the spatial-temporal incongruence between user travel demand and vehicle supply by reallocating idle vehicles \citep{mao2020dispatch}.
This can significantly increase vehicle utilization and reduce infrastructure capacities \citep{zhang2015exploring, zhang2017parking}.

As SAVs will revolutionize transportation systems, we urgently need to explore the day-to-day operation and mid- to long-term design of SAV systems.
When designing an SAV system, the key factors include long-term network design and parking infrastructure deployment planning, mid-term fleet sizing strategies, and short-term SAV dynamic routing and ridesharing operations.
Since long-term, mid-term, and short-term decision-making are interrelated, a unified framework that simultaneously tackles these planning, strategic, and operational problems must be essential \citep{li2021time, seo2021multi}.
Notably, SAV systems should be designed to explicitly capture trade-off relations among performance indexes associated with each level of the problem.
For example, larger fleet sizes help maintain a high service level for travelers, such as by reducing ride-share matching friction, however, they also require significant costs to purchase vehicles. 
Furthermore, operating a large fleet size without heavy congestion requires investment in road and parking capacity expansion.
Allocating parking slots in areas away from CBDs reduces infrastructure construction costs, yet at the expense of travelers' waiting time.
These examples illustrate that the designs of SAV systems are different depending on the strategic goals of systems’ operator \citep{dandl2021regulating}.

The uncertainty inherent in traffic demand emphasizes the importance of capturing the interrelations among planning, strategic, and operational problems.
To maintain the level of service to travelers, proactive reallocation of idle SAVs requires information regarding future traffic demand.
The historical request data can provide insight into statistical demand trends; however, future on-demand requests can only be predicted stochastically at best.
To immediately respond to stochastic on-demand demand, vehicles should be allocated close to the requested traveler.
Ensuring the feasibility of these operational processes requires significant investments in SAV and infrastructure resources.
Since investments in such supply resources are usually irreversible, the decision-making process is faced with trade-off relations between the service levels for uncertain future traffic demand and the large investment costs \citep{karoonsoontawong2007robust,ukkusuri2009multi}.

In addition to on-demand requests, pre-booking options may alleviate the impact of uncertainty on fleet size strategies and network design planning.
Travelers send information related to their trip request (e.g., departure time) in advance of daily SAV operation.
SAV operators make decisions on the number of SAVs to deploy in routing and ridesharing matching operations to meet the traffic demand of deterministic pre-booked and stochastic on-demand requests.
To provide an incentive for pre-booked requests, a fleet of vehicles dedicated to pre-booked travelers may be deployed.
This pre-booked option reduces the uncertainty faced in fleet sizing strategy decision-making and could have an indirect, though not direct, impact on network design planning.
\citet{yu2023coordinating} and \citet{chen2024real} demonstrated that integrating pre-booked requests into daily SAV operations can significantly increase the system performance.
Their work overlooks mid-term fleet sizing strategies and long-term network design, as well as trade-off relations with short-term operations.

This study focuses on SAV system design and operation problems with mixed traffic demand: on-demand and pre-booked requests.
The contributions of this paper can be summarized as follows:
\begin{itemize}
	\setlength{\itemsep}{0pt} \setlength{\parskip}{0pt}
	\item We formulate SAV system design and operation problems with pre-booked and on-demand requests as multi-stage stochastic linear programming (MSSLPs).
The MSSLPs jointly optimize macroscopic continuous variables regarding dynamic routing and ridesharing matching, fleet sizing, and network design.
It designs traffic capacities for roads and parking lots in the first stage and then determines the number of SAVs in the second stage. Pre-booked requests are known in the second phase. The routing and ridesharing matching of SAVs is then operated to satisfy stochastically on-demand requests along with the pre-booked requests.
To explicitly capture the dynamic operations and congestion propagation of SAV systems, SAV and traveler movements are expressed with a dynamic traffic assignment (DTA) framework as aggregate flows.
Our formulation can be considered as a stochastic extension of the linear programming for system optimum DTA problems for SAV systems \citep{levin2017congestion, li2021time, seo2021multi}.
    \item The MSSLPs can be decomposed into time-staged subproblems based on the Bellman Equations.
Leveraging the linearity of the MSSLPs, these subproblems can be rewritten as linear programming with piecewise linear value functions.
Furthermore, the linearity allows us to solve multi-objective optimization problems using the weighted-sum method.
Thanks to the linearity, a stochastic dual dynamic programming \citep{pereira1991multi} can be applied, which iteratively computes subproblems in bidirectional processes to obtain a solution with guaranteed optimality.
    \item Our numerical examples, one-dimensional urban networks, provide a sensitivity analysis of pre-booking options and dedicated vehicles, as well as an evaluation of the Pareto efficiency.
    We also confirm the performance of the SAV system with pre-booked requests by comparing the system only involving on-demand requests as a benchmark.
\end{itemize}

The remainder of this paper is organized as follows.
Section 2 introduces a SAV system with pre-booked and on-demand requests, and then formulates the corresponding MSSLP.
Section 3 explains SDDPs, typical solution procedures for MSSLPs, and then in Section 4, numerical examples provide the sensitivity analysis on pre-booking options and dedicated SAVs to serve only pre-booked travelers. 
This section also explains the Pareto-efficiency of the SAV system.
Finally, Section 5 summarizes the main achievements and suggestions for future work.

\section{Multi-stage Stochastic Linear Programming Formulation}

\subsection{Problem Description}
\label{sec:problem_description}
The SAV system considered in this study comprises five components: planning horizon, network, traveler, SAV, and a single SAV operator.
The system can be optimized when the total expected cost incurred during the planning horizon is minimized.
The following subsections describe the problem setting for each component.
We note that the basic structure of the problem other than stochastic demand is identical to that of \citet{seo2021multi}.

\subsubsection{Planning Horizon}
This study considers a finite and discrete planning horizon $\mathcal{T}_0=[0,...,T]$.
The length of the planning horizon $T$ is given.
The planning horizon $\mathcal{T}_0$ consists of an infrastructure construction stage $0$, a fleet size strategy stage $1$, and SAV operation stages $\mathcal{T}=[2,..,T]$.
Travelers send pre-booked requests in the fleet size strategy stage, 
whereas other travelers send on-demand requests at the beginning of their own desired departure time step within the SAV operation stages.
The length of the planning horizon $T$ is sufficiently large that both the pre-booked and on-demand requests can be satisfied within the planning horizon.

\begin{figure}[!ht]
  \centering
  \includegraphics[width=0.75\textwidth]{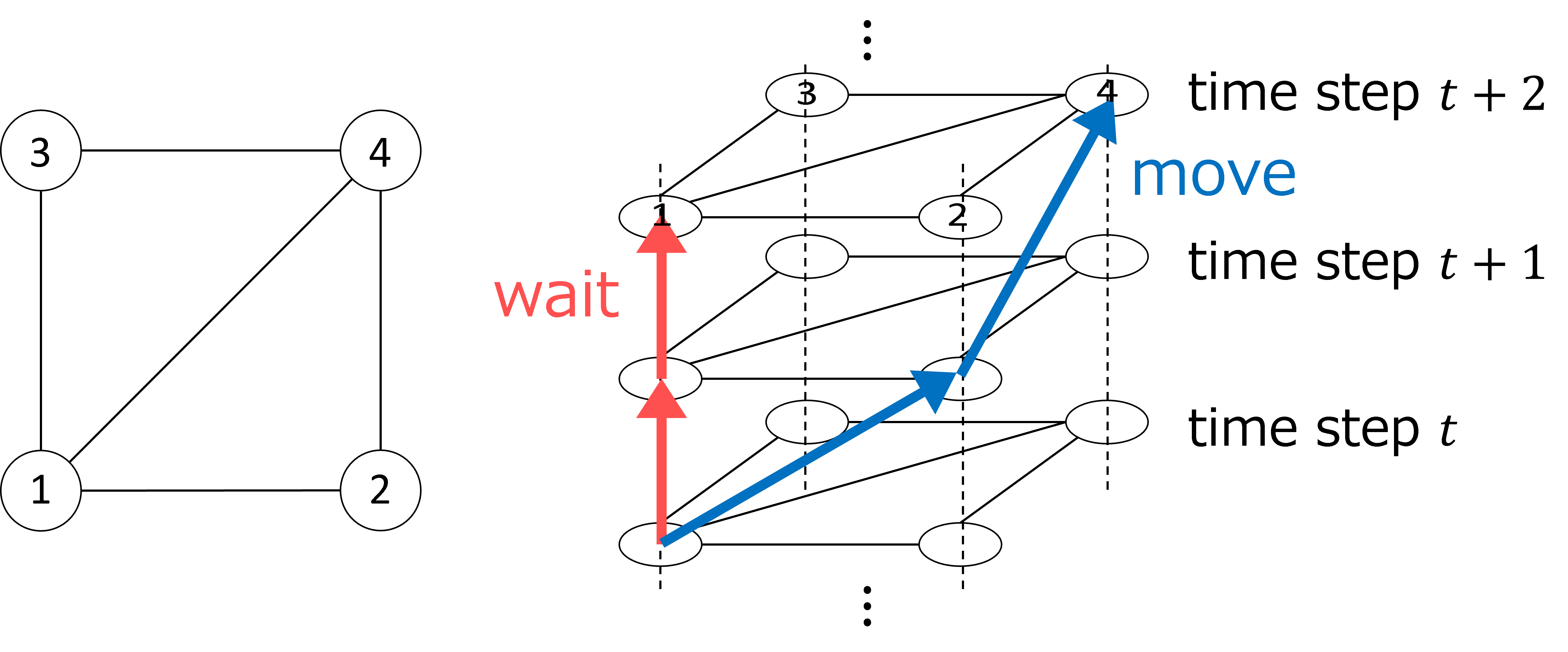}
  \caption{A Concept of Time-expanded Network} \label{fig:network}
\end{figure}

\subsubsection{Time-expanded Network}

We adopt a time-expanded network to describe the dynamic features of SAVs.
A time-expanded network is a road network representing the movement space with the addition of a time dimension.
Figure~\ref{fig:network} shows the concept of a time-expanded network.
For example, flows from node $1$ to node $2$ represent movement and travel time on link $(1,2)$, whereas flows from node $1$ to node $1$ represent waiting at node $1$.
We denote a set of nodes $\mathcal{N}$ in the time-expanded network and $\mathcal{L}$ be a set of links.
$\mathcal{I}_i$ and $\mathcal{O}_i$ are defined as a set of downstream and upstream nodes of node $i\in\mathcal{N}$, respectively.
A set of origin and destination nodes of traffic demand are $\mathcal{R}$ and $\mathcal{S}$, respectively.

All links have free-flow travel time and traffic capacity.
The free-flow Travel time of a link $ij\in\mathcal{L}$, $t_{ij}$, is constant, and the traffic capacity $\mu_{ij}$ is determined in the infrastructure construction stage $0$.
The traffic capacity of link $ii$, $\mu_{ii}$, represents the storage capacity of node $i\in\mathcal{N}$.
SAVs on node $i\in\mathcal{N}$ can flow into downstream node $j\in\mathcal{I}_i~(i\neq j)$ with the travel time $t_{ij}$ as long as its traffic capacity $\mu_{ij}$ is not exceeded.
SAVs that cannot flow into downstream node $j$ wait at node $i$ as long as its traffic capacity $\mu_{ii}$ is not exceeded.
This traffic flow model is equivalent to a point queue model with a finite queue length.

\subsubsection{SAV Operator}
The SAV operator at first determines traffic capacities $\mu_{ij}$ on each link $ij\in\mathcal{L}$ and $\mu_{ii}$ on each node $i\in\mathcal{N}$ in the infrastructure construction stage $0$.
He/she then deploys SAVs while observing pre-booked requests from travelers.
Subsequently, the deployed SAVs run and stop to pick-up and drop-off travelers on the time-expanded network, following the optimized routes and schedules while observing on-demand requests.
The SAV operator can adapt fleet sizing strategies and SAV operations to minimize the total expected cost according to the observed trip requests.
Whereas, the network design planning is forced to be determined under unknown and stochastic pre-booked and on-demand requests.
All decision variables can be optimized when the total expected cost, including the total travel distance and the total number of SAVs, the infrastructure expansion cost, and the total travel time of travelers.

\subsubsection{Traveler}
Travelers send pre-booked or on-demand requests to the SAV operator.
Trip requests are composed of specific origin $r\in\mathcal{R}$, destination $s\in\mathcal{S}$, earliest departure time $k\in\mathcal{K}\subset\mathcal{T}$, and latest arrival time $T^k$.
 $\mathcal{K}$ is a set of time steps corresponding to travelers' departure time,
and $T^k$ is a time step corresponding to the latest arrival time of travelers with departure time step $k$.
We assume that the time value of travelers is homogeneous and that the probability distributions of both trip requests are given and stage-wise independent.

They travel only by SAVs according to the optimized route suggested by the SAV operator and eventually arrive at their destinations.
The number of travelers that a SAV can carry is constrained by a pre-determined vehicle carrying capacity.
If Travelers are unable to ride SAVs, then they wait at nodes.
The length of the traveler's queue is not restricted.

\subsubsection{SAV}
The dynamics of SAV flow is expressed using a point queue model with a limited queue length.
SAVs move on links at the free flow speed and stop at nodes when parked or in traffic congestion.
The congestion occurs when the traffic volume reaches the traffic capacity of the link or the waiting capacity of the node. 
SAVs follow the optimized itinerary and route directed by the SAV operator for all choices, such as starting and stopping.
SAVs always move on roads or park in parking lots during the planning horizon.
SAVs can pick up passengers as long as the pre-determined vehicle carrying capacity is not reached.
This study considers two types of SAVs: dedicated SAVs delivering only pre-booked travelers and non-dedicated SAVs accommodating both pre-booked and on-demand travelers.

\begin{table*}
	\caption{List of variable notation}
	\label{variable}       
	\scalebox{0.95}[0.95]{
		\begin{tabular}{ll}
			\hline\noalign{\smallskip}
			notation & definition \\
			\noalign{\smallskip}\hline\noalign{\smallskip}
			$x_{ij}^{t}$ & flow of dedicated SAVs that start traveling link $ij\in\mathcal{L}$ on time step $t\in\mathcal{T}$\\
			$\hat{x}_{ij}^{t}$ & flow of non-dediicated SAVs that start traveling link $ij\in\mathcal{L}$ on time step $t\in\mathcal{T}$\\
			$y_{s,ij}^{k,t}$ & flow of pre-booked travelers who start traveling link $ij\in\mathcal{L}$ on time step $t\in\mathcal{T}^k$, destination node $s\in\mathcal{S}$, and\\
			&departure time step $k\in\mathcal{K}$ \\
			$\mathcal{Y}_{s,ij}^{k,t}$ & flow of pre-booked travelers who start traveling link $ij\in\mathcal{L}$ in dedicated SAVs, on time step $t\in\mathcal{T}^k$, \\
			&destination node $s\in\mathcal{S}$, and departure time step $k\in\mathcal{K}$ \\
			$\hat{y}_{s,ij}^{k,t}$ & flow of on-demand travelers who start traveling link $ij\in\mathcal{L}$ on time step $t\in\mathcal{T}^k$, destination node $s\in\mathcal{S}$, and \\
			&departure time step $k\in\mathcal{K}$\\
			$Q_s^{k,t}$& Cumulative number of pre-booked traveler arrivals on time step $t\in\mathcal{T}^k$, with destination node $s\in\mathcal{S}$\\
			&and departure time step $k\in\mathcal{K}$\\
			$\hat{Q}_s^{k,t}$& Cumulative number of on-demand traveler arrivals on time step $t\in\mathcal{T}^k$, with destination node $s\in\mathcal{S}$\\
			& and departure time step $k\in\mathcal{K}$\\
			$\mu_{ij}$ & traffic capacity of link $ij\in\mathcal{L}$\\
			$T^t$ & instantaneous total travel time of travelers on time step $t\in\mathcal{T}$ (including waiting time on nodes)\\
			$D^t$ & instantaneous total distance traveled by SAVs on time step $t\in\mathcal{T}$\\
			$N$ & total number of SAVs\\
			$C$ & total cost of infrastructure construction\\
			\noalign{\smallskip}\hline
	\end{tabular}}
\end{table*}

\begin{table*}
	\caption{List of parameter notation}
	\centering
	\label{parameter}       
	\scalebox{0.85}[0.85]{
		\begin{tabular}{ll}
			\hline\noalign{\smallskip}
			notation & definition \\
			\noalign{\smallskip}\hline\noalign{\smallskip}
			$t_{ij}$ & free-flow travel time of link $ij\in\mathcal{L}$ if $i\neq{j}$\\
			$t_{ii}$ & waiting time at node $i\in\mathcal{N}$ for one time step (i.e., equal to the time step width by the definition)\\
			$d_{ij}$ & length of link $ij\in\mathcal{L}$\\
			$c_{ij}$ & unit cost of expanding traffic capacity of link $ij\in\mathcal{L}$ if $i\neq{j}$\\
			$c_{i}$ & unit cost of expanding storage capacity of node $i\in\mathcal{N}$ for SAVs \\
			$\rho$ & carrying capacity of an SAV\\
			$\mu_{ij}^\text{max}$, $\mu_{ij}^\text{min}$ & maximum and minimum allowable value of $\mu_{ij}$, respectively\\
			$\xi_{rs}^{k}$ & time-dependent pre-booked requests of travelers with origin $r\in\mathcal{R}$, destination $s\in\mathcal{S}$, and departure time step $k\in\mathcal{K}$\\
			$\hat{\xi}_{rs}^{k}$ & time-dependent on-demand requests of travelers with origin $r\in\mathcal{R}$, destination $s\in\mathcal{S}$, and departure time step $k\in\mathcal{K}$ \\
			$T^k$ & latest arrival time of traveler with departure time step $k\in\mathcal{K}$\\
			$T$ & final time step\\
			\noalign{\smallskip}\hline
	\end{tabular}}
\end{table*}

\subsection{Formulation}
This section formulates the SAV system design and operation problem with pre-booked and on-demand requests according to the problem description described in the previous subsection.
The definitions of variables and parameters are summarized in Tables \ref{variable} and \ref{parameter}, respectively.

\begin{figure}[!t]
	\centering
	\includegraphics[width=0.9\textwidth]{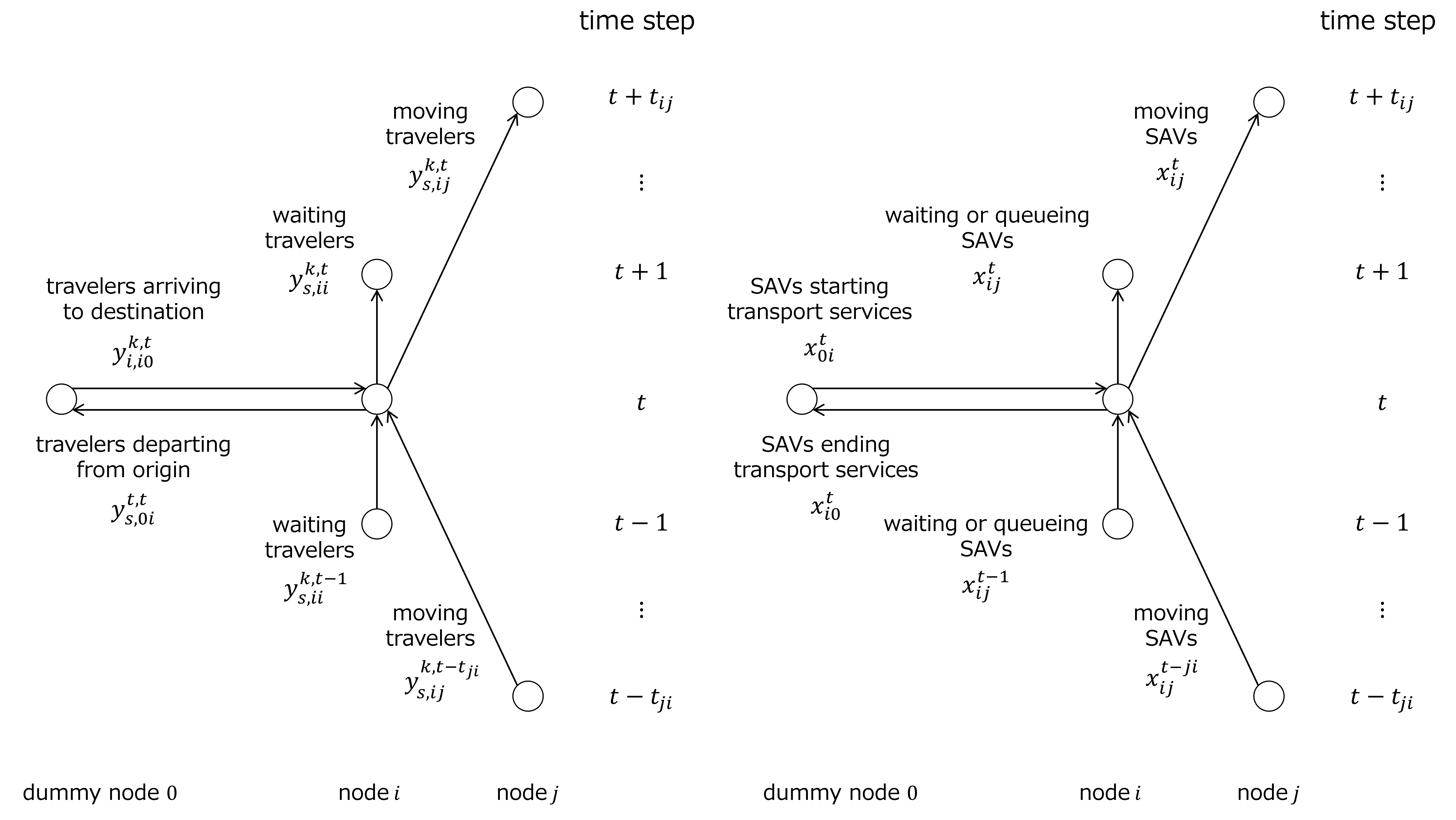}
	\caption{Network Flow} \label{fig:flow}
\end{figure}

\subsubsection{Network Flow}
SAV operations at time step $t\in\mathcal{T}$ are defined by SAV and traveler flow on the time-expanded network.
The dynamics of SAV flow is governed by traffic capacity and trip time and can be represented by a point queue model with a finite queue length.
That of traveler flow is restricted by the vehicle carrying capacity.
The flows must also satisfy conservation constraints at any node shown in Figure \ref{fig:flow}.
Denoting flow of dedicated and non-dedicated SAVs that travel link $ij$ on time step $t$ as $x_{ij}^t$ and $\hat{x}_{ij}^t$, respectively, and flow of pre-booked and on-demand travelers who start traveling link $ij$ on time step $t$, destination node $s$, and departure time step $k$ as $y_{s,ij}^{k,t}$ and $\hat{y}_{s,ij}^{k,t}$, respectively, flow constraints that should be ensured on a time-expanded network can be expressed as follows:
\begin{flalign}
	\nonumber
	&\underline{\rm Flow~Conservation~Constraints}&\\
	\label{eq:State_Equation_x}
	&\sum_{j\in\mathcal{O}_i}x_{ji}^{t-t_{ji}} +x_{0i}^{t-1}=\sum_{j\in\mathcal{I}_i}x_{ij}^{t}+x_{i0}^{t}&\forall i,t\in\mathcal{T},\\
	\label{eq:State_Equation_hatx}
	&\sum_{j\in\mathcal{O}_i}\hat{x}_{ji}^{t-t_{ji}} +\hat{x}_{s,0i}^{t-1}=\sum_{j\in\mathcal{I}_i}\hat{x}_{ij}^{t}+\hat{x}_{i0}^{t}&\forall i,t\in\mathcal{T},\\
	\label{eq:State_Equation_y}
	&\sum_{j\in\mathcal{O}_i}y_{s,ji}^{k,t-t_{ji}} +y_{s,0i}^{k,t}=\sum_{j\in\mathcal{I}_i}y_{s,ij}^{k,t}+y_{s,i0}^{k,t}&\forall i,k,s,t\in\mathcal{T}^k,\\
	\label{eq:State_Equation_haty}
	&\sum_{j\in\mathcal{O}_i}\hat{y}_{s,ji}^{k,t-t_{ji}} +\hat{y}_{s,0i}^{k,t}=\sum_{j\in\mathcal{I}_i}\hat{y}_{s,ij}^{k,t}+\hat{y}_{s,i0}^{k,t}&\forall i,k,s,t\in\mathcal{T}^k,\\
	\nonumber
	&\underline{\rm Demand~Constraints}&\\
	\label{eq:State_Equation_q}
	&Q_{s}^{k,t}=Q_{s}^{k,t-1}+y_{s,s0}^{k,t}&\forall k,s,t\in\mathcal{T}^k,\\
	\label{eq:State_Equation_q_hat}
	&\hat{Q}_{s}^{k,t}=\hat{Q}_{s}^{k,t-1}+\hat{y}_{s,s0}^{k,t}&\forall k,s,t\in\mathcal{T}^k,
 \end{flalign}
\vspace{-7ex}
\begin{flalign}
	\label{eq:Origin_Demand}
	&y_{s,0r}^{k,t}=\begin{cases}
		\xi_{rs}^t & {\rm if}~t = k \\
		0 & {\rm otherwise}
	\end{cases}&\forall k,rs,t\in\mathcal{T}^k,\\
	\label{eq:Origin_Demand_hat}
	&\hat{y}_{s,0r}^{k,t}=\begin{cases}
		\hat{\xi}_{rs}^t & {\rm if}~t = k \\
		0 & {\rm otherwise}
	\end{cases}&\forall k,rs,t\in\mathcal{T}^k,\\
	\label{eq:Destination_Demand}
	&Q_{s}^{k,t}=\sum_{r\in\mathcal{R}}\xi_{rs}^k&\forall k,s,t\in[T^k,T],\\
	\label{eq:Destination_Demand_hat}
	&\hat{Q}_{s}^{k,t}=\sum_{r\in\mathcal{R}}\hat{\xi}_{rs}^k&\forall k,s,t\in[T^k,T],\\
	\nonumber
	&\underline{\rm Supply~Constraints}&\\
	\label{eq:FleetSize}
	&N=\sum_i \{x_{0i}^1+\hat{x}_{0i}^1\},&\\
	\label{eq:SAV_supply}
	&x_{0i}^t=0&\forall i,t\in\mathcal{T}, \\
	\label{eq:SAV_supply0}
	&\hat{x}_{0i}^t=0&\forall i,t\in\mathcal{T}, \\
	\label{eq:SAV_return}
	&x_{i0}^t=0&\forall i,t\neq T, \\
	\label{eq:SAV_return0}
	&\hat{x}_{i0}^t=0&\forall i,t\neq T, \\
	\nonumber
	&\underline{\rm Capacity~Constraints}&\\
	\label{eq:Road_Capacity_Constraint}
	&x_{ij}^{t}+\hat{x}_{ij}^{t}\le\mu^t_{ij}&\forall ij,t\in\mathcal{T},\\
	\label{eq:State_Equation_u}
	&\mu_{ij}^{t}-\mu_{ij}^{t-1}=0&\forall ij,t\in\mathcal{T}_0,\\
	\label{eq:min_max_u}
	&\mu_{ij}^{\rm min}\le \mu_{ij}^0 \le \mu_{ij}^{\rm max}&\forall ij,\\
	\label{eq:Vehicle_Capacity_Constraint}
	&\sum_{k,s}\mathcal{Y}_{s,ij}^{k,t}\le \rho x^t_{ij}&\forall ij,t\in\mathcal{T},\\
	\label{eq:Vehicle_Capacity_Constraint_hat}
	&\sum_{k,s}(y_{s,ij}^{k,t} - \mathcal{Y}_{s,ij}^{k,t}) + \sum_{k,s}\hat{y}_{s,ij}^{k,t}\le \rho \hat{x}^t_{ij}&\forall ij,t\in\mathcal{T},\\
	\label{eq:Passenger_Conservation}
	&0\le \mathcal{Y}_{s,ij}^{k,t}\le y_{s,ij}^{k,t}&\forall ij,k,s,t\in\mathcal{T}^k,\\
	&\underline{\rm Non\mathchar`-negative~Constraints}&\\
	\label{eq:Nonnegative_x}
	&x_{ij}^t\ge0&\forall ij, t\in\mathcal{T},\\
	\label{eq:Nonnegative_xhat}
	&\hat{x}_{ij}^t\ge0&\forall ij, t\in\mathcal{T},\\
	\label{eq:Nonnegative_y}
	&y_{s,ij}^{k,t}\ge0&\forall ij,k,s,t\in\mathcal{T}^k,\\
	\label{eq:Nonnegative_yhat}
	&\hat{y}_{s,ij}^{k,t}\ge0&\forall ij,k,s,t\in\mathcal{T}^k,\\
	\label{eq:Nonnegative_x0}
	&x_{0i}^1\ge0&\forall i,\\
	\label{eq:Nonnegative_xhat0}
	&\hat{x}_{0i}^1\ge0&\forall i,\\
	\label{eq:Nonnegative_y0}
	&y_{s,s0}^{k,t}\ge0&\forall k,s,t\in\mathcal{T}^k,\\
	\label{eq:Nonnegative_yhat0}
	&\hat{y}_{s,s0}^{k,t}\ge0&\forall k,s,t\in\mathcal{T}^k,\\
	\label{eq:Nonnegative_Q}
	&{Q}_{s}^{k,t}\ge0&\forall k,s,t\in\mathcal{T}^k,\\
	\label{eq:Nonnegative_Qhat}
	&\hat{Q}_{s}^{k,t}\ge0&\forall k,s,t\in\mathcal{T}^k.
\end{flalign}

Eqs.(\ref{eq:State_Equation_x})--(\ref{eq:State_Equation_haty}) are node conversation constraints for SAV and traveler flows.
The left and right sides in each equation describe total inflow and outflow, respectively.
Node $0$ is a dummy node representing departures (incoming into the network) and arrivals (outgoing from the network).
Eqs.(\ref{eq:State_Equation_q})--(\ref{eq:Destination_Demand_hat}) represent traffic demand constraints.
Eqs.(\ref{eq:State_Equation_q}) and (\ref{eq:State_Equation_q_hat}) are transition equations for cumulative arrivals of pre-booked and on-demand travelers, respectively.
Eqs.(\ref{eq:Origin_Demand}) and (\ref{eq:Origin_Demand_hat}) explain demand generation constraints for pre-booked and on-demand travlers, respectively.
Eqs.(\ref{eq:Destination_Demand}) and (\ref{eq:Destination_Demand_hat}) are demand attraction constraints for pre-booked and on-demand travlers, respectively.
Eqs.(\ref{eq:FleetSize})--(\ref{eq:SAV_return0}) represent traffic supply constraints.
Eqs.(\ref{eq:FleetSize})--(\ref{eq:SAV_supply0}) ensure that SAVs are deployed only in the fleet sizing strategy stage.
Eqs.(\ref{eq:SAV_return}) and (\ref{eq:SAV_return0}) ensure that SAVs exist on the time-extended network in the planning horizon.
Eqs. (\ref{eq:Road_Capacity_Constraint})--(\ref{eq:Passenger_Conservation}) is associated with capacity constraints.
Eq.(\ref{eq:Road_Capacity_Constraint}) means traffic capacity constraints.
Eq.(\ref{eq:State_Equation_u}) is the proxy equation that transfers the infrastructure construction decision to the next stage.
Eq.(\ref{eq:min_max_u}) describes the maximum and minimum allowable of traffic capacity.
Eqs.(\ref{eq:Vehicle_Capacity_Constraint}) and (\ref{eq:Vehicle_Capacity_Constraint_hat}) explain vehicle carrying capacity constraints for dedicated and non-dedicated SAVs, respectively.
Eq.(\ref{eq:Passenger_Conservation}) ensures that the maximum number of pre-booked travelers available to board the dedicated SAVs is equal to total pre-booked travelers.
Eqs.(\ref{eq:Nonnegative_x})--(\ref{eq:Nonnegative_Qhat}) represent non-negative constraints on all variables.

\subsubsection{Objective Function}

SAV operator aims to minimize the total expected cost incurred during the planning horizon through infrastructure design planning, fleet sizing strategy, and SAV operations.
The total expected cost is the time-integrals of the instantaneous costs $F^t$ incurred in the SAV system on time step $t\in\mathcal{T}_0$.
The instantaneous costs $F^t$ are composed of the total travel time, $T^t$, the total distance traveled by SAVs, $D^t$, the total number of SAVs, $N$, and the total infrastructure cost $C$, and given by
\begin{flalign}
	\label{eq:TravelTime}
	&T^t=\sum_{ij,k,s}\{t_{ij}y_{s,ij}^{k,t}+t_{ij}\hat{y}_{s,ij}^{k,t}\}&\forall t\in\mathcal{T},\\
	\label{eq:Distance}
	&D^t=\sum_{ij, i\neq j}\{d_{ij}x_{ij}^t + d_{ij}\hat{x}_{ij}^t\}&\forall t\in\mathcal{T},\\
	&N=\sum_i \{x_{0i}^1+\hat{x}_{0i}^1\},\tag{\ref{eq:FleetSize}}&\\
	\label{eq:InfraCost}
	&C=\sum_{ij, i\neq j}c_{ij}\mu^0_{ij} + \sum_i c_i\mu^0_{ii}.&
\end{flalign}

This study presents a multi-objective optimization problem that jointly minimizes the total travel time of travelers, the total distance traveled by SAVs, the total number of SAVs, and the total infrastructure cost.
Solving a multi-objective optimization problem involves deriving its Pareto frontier---a set of Pareto-efficient solutions \citep{ehrgott2005multicriteria}.
Since system constraints as shown in Eqs.(\ref{eq:State_Equation_x})--(\ref{eq:InfraCost}) are linear,  
the Pareto-efficient solutions can be calculated by solving single objective optimization problems weighted with appropriate non-negative weights using the weighted sum method.
Let the weights of each performance index be $\alpha_T$, $\alpha_D$, $\alpha_N$, and $\alpha_C$, objective function can be obtain as follows:
\begin{flalign}
	\label{eq:ObjectiveFunction}
	&F^t=\begin{cases}
		\alpha_C C & {\rm if}~~t = 0\\
		\alpha_N N & {\rm if}~~t = 1\\
		\alpha_T T^t + \alpha_D D^t & {\rm otherwise} 
	\end{cases}&\forall t\in\mathcal{T}_0.&
\end{flalign}

\subsubsection{Multi-stage Stochastic Linear Program}

This subsection presents the multi-stage stochastic linear program (MSSLP) to find optimal solutions 
subject to constraints based on the problem description and network flow described in the previous section.
The optimal solutions can be decomposed into infrastructure planning, fleet sizing, and SAV operations.
Optimal infrastructure planning minimizes construction costs and expected future costs within the planning period, under unknown and stochastic traffic demand.
Optimal fleet size can be determined while observing pre-booked requests.
Optimal SAV operation should be flexible in response to stochastic on-demand demand.
The system optimal state can be derived by the following multi-stage stochastic programming problem:
\begin{flalign}
	\nonumber
	&\underline{\rm [SAV\mathchar`-MSSLP]}&\\
	\label{eq:Original_Objective_Function}
	&\min_{\bm{z}^0\in\mathcal{X}^0}	F^0(\bm{z}^0)
	+ \mathbb{E}^{{1}}\left[\min_{\bm{z}^1\in\mathcal{X}^1(\bm{z}^0,\bm{\xi}^1)}	F^1(\bm{z}^1,\bm{\xi}^1)
	+\mathbb{E}^{{2}}\left[\cdots+
	\mathbb{E}^{{T}}\left[\min_{\bm{z}^T\in\mathcal{X}^T(\bm{z}^{T-1},\bm{\xi}^{T})}	F^T(\bm{z}^T,\bm{\xi}^T)
	\right]	\right]	\right],&
\end{flalign}
where $\mathbb{E}^t[\cdot]$ denotes the conditional expectation given the information up to time step $t$,
$\bm{z}^t$ is a vector of the decision variable at time step $t$, 
$\bm{\xi}^t$ represents a vector of traffic demand, including pre-booked and on-demand requests, observed at time step $t$.
$\mathcal{X}^0$ is the feasible region of decision variables before observing any traffic demand,
$\mathcal{X}^t(\bm{z}^{t-1},\bm{\xi}^t)$ is the feasible region given the past decision $\bm{z}^{t-1}$ and traffic demand $\bm{\xi}^t$.
These feasible regions are given by the following linear constraints:
\begin{flalign}
	&\mathcal{X}^0=\Bigl\{\bm{\mu}^0\bigl|{\rm Eqs.}(\ref{eq:min_max_u}),(\ref{eq:InfraCost}) {\rm~and~} (\ref{eq:ObjectiveFunction})\Bigr\},&\\
	&\mathcal{X}^1(\bm{z}^{0},\bm{\xi}^1)=\Bigl\{\bm{z}^1=(\bm{x}^1,\bm{\mu}^1)\bigl|~{\rm Eqs.(\ref{eq:FleetSize}),(\ref{eq:State_Equation_u}),(\ref{eq:Nonnegative_x0}),(\ref{eq:Nonnegative_xhat0}), and~  (\ref{eq:ObjectiveFunction})}\Bigr\},&\\
	&\mathcal{X}^t(\bm{z}^{t-1},\bm{\xi}^t)=\Bigl\{\bm{z}^t=(\bm{x}^t,\bm{y}^t,\bm{Q}^t,\bm{\mu}^t)\bigl|~{\rm Eqs.(\ref{eq:State_Equation_x})\mathchar`-(\ref{eq:ObjectiveFunction})}\Bigr\}&\forall t\in\mathcal{T},
\end{flalign}
where $\bm{\mu}$, $\bm{x}$, $\bm{y}$, and $\bm{Q}$ are vectors of traffic capacities, SAV flows, traveler flows, and cumulative arrivals, respectively.

\section{Stochastic Dual Dynamic Programming}

Stochastic dual dynamic programming (SDDP) is a sampling-based algorithm for solving the MSSLPs under the assumption of stage-wise independence.
Since the classical SDDP had been developed by \citet{pereira1991multi}, 
SDDP has been applied in a variety of fields, including energy system scheduling (\citep{homem2011sampling,shapiro2013risk,maceiral2018twenty}), transportation planning (\citep{fhoula2013stochastic,angun2015stochastic,delgado2019multistage}), and financial planning (\citep{kozmik2015evaluating,homem2016risk,lee2023large}).
This section outlines a generic SDDP algorithm for solving the [SAV-MSSLP].
More comprehensive discussions can be found in \citet{philpott2008convergence} and \citet{shapiro2011analysis}.

\subsection{Generic Formulation of MSSLPs}

Let us consider the following a generic MSSLP:
\begin{subequations}
	\label{eq:MSSP}
	\begin{flalign}
		&\min_{\bm{x}_1\ge\bm{0}}\bm{c}_1^{\top}\bm{x}_1+\mathbb{E}_{P_2|\bm{\xi}_1}\left[\min_{\bm{x}_2\ge\bm{0}}\bm{c}_2^{\top}\bm{x}_2+\mathbb{E}_{P_3|\bm{\xi}_2}\left[\cdots+\mathbb{E}_{P_{T}|\bm{\xi}_{T-1}}\left[\min_{\bm{x}_T\ge\bm{0}}\bm{c}_T^{\top}\bm{x}_T\right]\right]\right],&\label{eq:MSSP_objective}\\
		&\bm{A}_t\bm{x}_t+\bm{B}_t\bm{x}_{t-1}=\bm{b}_t,~{\rm where~}\bm{x}_0 {\rm~and~} \bm{\xi}_1{\rm~are~given,}&\label{eq:MSSP_transition}\forall t\in\mathcal{T}.
	\end{flalign}
\end{subequations}
where the vectors $\bm{b}_t,\bm{c}_t$ and the matrices $\bm{A}_t,\bm{B}_t$ are random variables forming the stochastic process $\bm{\xi}_t=(\bm{b}_t,\bm{c}_t,\bm{A}_t,\bm{B}_t)$.
$\mathbb{E}_{P_t|\bm{\xi}_{[t-1]}}[\cdot]$ denotes the conditional expectation given the past observations $\bm{\xi}_{[t-1]}$ under the joint probability distribution $P_t$.
We describe the SDDP algorithm to solve the MSSLP (\ref{eq:MSSP}).

The MSSLP (\ref{eq:MSSP}) can be decomposed into the following subproblems based on the Bellman optimality principle:
\begin{subequations}
	\label{eq:MSSP_DP}
	\begin{flalign}
		&V_{t}(\bm{x}_{t-1},\tilde{\bm{\xi}}_t)=\min_{\bm{x}_t\in\bm{\mathcal{X}}_t(\bm{x}_{t-1},\tilde{\bm{\xi}}_t)}\left[\bm{c}_t^{\top}\bm{x}_t+\mathcal{V}_{t+1}(\bm{x}_t)\right]&\forall t,\label{eq:MSSP_DP_a}\\
		&\mathcal{V}_{t+1}(\bm{x}_{t})=\mathbb{E}_{P_{t+1}|\tilde{\bm{\xi}}_{t}}\left[V_{t+1}(\bm{x}_{t},{\bm{\xi}}_{t+1})\right] {\rm~where~}\mathcal{V}_{T+1}(\cdot)=0&\forall t.\label{eq:MSSP_DP_b}
	\end{flalign}
\end{subequations}
where $V_{t}(\bm{x}_{t-1},\bm{\tilde{\xi}}_{t})$ is the future cost at time step $t$ and $\mathcal{V}_{t}(\bm{x}_{t-1})$ is the expected future cost.
$\tilde{\bm{\xi}_t}$ is the realization of the random variable.
We should note that if $\bm{\xi}_{t}$ is continuous, the computation of Eq.(\ref{eq:MSSP_DP_b}) requires multi-dimensional integral calculations.
A promising approach is to solve the problem with a Sample Average Approximation (SSA) constructed from random sampling from the true probability distribution.
According to \citet{shapiro2011analysis}, if SDDP iterations and the number of samples are sufficient, the optimal solutions of the SSA problem and the original problem coincide.
A random sample of the stochastic process $\bm{\xi}_{t}$ is expressed as $\bm{\xi}_t=\{\bm{\tilde{\xi}}_t^1,.... \bm{\tilde{\xi}}_t^{N_t}\}$, 
where $N_t$ is the number of samples at time step $t$.
Replacing the true probability distribution by a probability distribution constructed from random samples, each with probability $1/N_t$, Eq.(\ref{eq:MSSP_DP}) can be rewritten as follows:
\begin{subequations}
    \label{eq:MSSP_DP_SSA}
    \begin{flalign}
		&{V}_{t}(\bm{x}_{t-1},\bm{\tilde{\xi}}^n_{t}) =\min_{\bm{x}_t\in\bm{\mathcal{X}}_t(\bm{x}_{t-1},\tilde{\bm{\xi}}^n_t)}\left[(\bm{c}_t^{n})^{\top}\bm{x}_t+\mathscr{V}_{t+1}(\bm{x}_t)\right]&\label{eq:MSSP_DP_SSA_a}\forall t,n=1,...,N_t,\\
		&\mathscr{V}_{t}(\bm{x}_{t-1})=\frac{1}{N_t}\sum_{n=1}^{N_t}\left[{V}_{t}(\bm{x}_{t-1},\bm{\tilde{\xi}}_{t}^n)\right],{\rm~where~}\mathscr{V}_{T+1}(\bm{x}_{T})=0&\forall t. \label{eq:MSSP_DP_SSA_b}
    \end{flalign}
\end{subequations}

\subsection{Generic Algorithm}

SDDP algorithm consists of a forward step to evaluate policies and a backward step to improve policy.
The backward step updates information about expected future costs by going back from the final step to the initial step and solving a subproblem at each time step.
This information is obtained as a piecewise linear approximation of the expected future cost in the form of a Benders cut.
Under loose conditions (e.g., convexity of the instantaneous cost), the piecewise linear function defined by a set of the Benders cuts provides an outer approximation (i.e., lower bound) of the expected future cost.
This avoids the curse of dimensionality arising from the discretization of variables.
The forward step evaluates the upper bound of the expected future cost by generating realizations (sample paths) of the stochastic process and computing the optimal solution corresponding to the realizations from the initial to the final step.
The Benders cut set obtained in the backward step is used as the expected future cost at each time step.
The calculation procedure is shown in Algorithm \ref{alg:overview}.

The iterative computation converges with probability 1 if the following assumptions are satisfied \citep{philpott2008convergence,shapiro2011analysis,dowson2021sddp}:
\begin{enumerate}
	\setlength{\itemsep}{0pt} \setlength{\parskip}{0pt}
	\item Given a realization of the stochastic process, the instantaneous cost function at each time step is convex, the transition equation is linear, and the feasible region is a nonempty bounded convex set,
	\item Stochastic process is stage-wise independent,
	\item Sample space is finite,
	\item Relatively complete recourse is satisfied, and
	\item The planning period is a finite or a time-discounted infinite.
\end{enumerate}

\begin{figure}[!t]
	\begin{algorithm}[H]
		\caption{An overview of the SDDP algorithm}
		\label{alg:overview}
		\begin{algorithmic}[]
			\setlength{\itemsep}{0pt} \setlength{\parskip}{0pt}
			\REQUIRE $\varepsilon > 0~({\rm convergence~condition}),~\{\tilde{\bm{\xi}}^1_t,...,\tilde{\bm{\xi}}^N_t\}_{t=2,...,T}$
			\STATE Initialize: iteration$\leftarrow0, \underline{\Psi}=-\infty~({\rm lower~bound}), \overline{\Psi}=\infty~({\rm upper~bound}), \{{\bm{\mathcal{K}}}_t=\emptyset\}_{t=1,...,T-1}~({\rm cut~index~sets})$
			\WHILE{$\overline{\Psi} - \underline{\Psi}>\varepsilon$}
			\STATE Perform the forward step: Update the upper bound $\overline{\Psi}$
			\STATE Perform the backward step: Update the lower bound $\underline{\Psi}$
			\STATE Update the iternation number: iteration$\leftarrow$iteration$+1$
			\ENDWHILE
		\end{algorithmic}
	\end{algorithm}
\end{figure}

\subsubsection{Backward Step}
\label{backward}

The backward step solves the subproblem backward using the $M$ trial solutions $\{\overline{\bm{x}}_t^m\}_{t=0,...,T-1,m=1,...,M}$ computed in the forward step.
In the following, we outline the computational procedure for the $m$-th trial solution to avoid notational complexity.

Let us consider a subproblem at the final stage $t=T$.
Given a trial solution $\overline{\bm{x}}_{T-1}$, Eq.(\ref{eq:MSSP_DP_SSA}) can be rewritten as follows:
\begin{subequations}
	\label{eq:backwardT}
	\begin{flalign}
		&{V}_{T}(\overline{\bm{x}}_{T-1},\bm{\tilde{\xi}}^n_{T}) = \min_{\bm{x}_T\ge0, \bm{A}^n_T\bm{x}_T=-\bm{B}^n_T\overline{\bm{x}}_{T-1}+\bm{b}^n_T}(\bm{c}_T^{n})^{\top}\bm{x}_T &\forall n=1,...,N_T,\label{eq:HJB_backwardT}\\
		&\mathscr{V}_{T}(\overline{\bm{x}}_{T-1})=\frac{1}{N_T}\sum_{n=1}^{N_T}\left[{V}_{T}(\overline{\bm{x}}_{T-1},\bm{\tilde{\xi}}_{T}^n)\right].\label{eq:ExpQ_backwardT}&
	\end{flalign}
\end{subequations}
Eq.(\ref{eq:HJB_backwardT}) is assumed to be well-defined for almost all $\bm{\xi}_{T}$ realizations of the optimal solution ${V}_{T}(\overline{\bm{x}}_{T-1},\bm{\tilde{\xi}}^n_{T})$ (i.e., satisfying a relatively complete recourse).
Since $\mathscr{V}_{T}(\overline{\bm{x}}_{T-1})$ is a convex function (actually a piecewise linear function), the Benders cut satisfies the following inequality:
\begin{flalign}
	\label{eq:Benders_inequality_multi}
	&{V}_{T}({\bm{x}}_{T-1},\bm{\tilde{\xi}}^n_{T})\ge{V}_{T}(\overline{\bm{x}}_{T-1},\bm{\tilde{\xi}}^n_{T})+(\bm{g}_T^n)^{\top}(\bm{x}_{T-1}-\overline{\bm{x}}_{T-1})&\forall n=1,...,N_T,\\
	\label{eq:Benders_inequality_single}
	&\mathscr{V}_{T}({\bm{x}}_{T-1})\ge\mathscr{V}_{T}(\overline{\bm{x}}_{T-1})+\bm{g}_T^{\top}(\bm{x}_{T-1}-\overline{\bm{x}}_{T-1}).&
\end{flalign}
where $\bm{g}_T$ and $\bm{g}_T^n$ are subgradients of $\mathscr{V}_{T}({\bm{x}}_{T-1})$ and ${V}_{T}({\bm{x}}_{T-1},\bm{\tilde{\xi}}^n_{T})$ at $\overline{\bm{x}}_{T-1}$, respectively,
using the dual variable vector $\bm{\pi}_T^n$ corresponding to the transition equation, given by
\begin{flalign}
	&\bm{g}_T^n=-(\bm{B}_T^{n})^{\top}\bm{\pi}_T^n&\forall n=1,...,N_T,\\
	&\bm{g}_T=\frac{1}{N_T}\sum_{n=1}^{N_T}\bm{g}_T^n.&
\end{flalign}

From Eqs.(\ref{eq:Benders_inequality_multi})(\ref{eq:Benders_inequality_single}), the Benders cut $\ell({\bm{x}}_{T-1})$ is given by
\begin{flalign}
    &\ell({\bm{x}}_{T-1}):=\bm{g}_T^{\top}\bm{x}_{T-1}+\beta_{T},&\\
    &\beta_{T}=\mathscr{V}_{T}(\overline{\bm{x}}_{T-1})-\bm{g}_T^{\top}\overline{\bm{x}}_{T-1}.&
\end{flalign}
Because of the stage-wise independence of stochastic processes, this Benders cut can be shared by all problems in the $(T-1)$ step.

Then, consider a subproblem in the $t=T-1$ stage.
The index set of Benders cuts formed by the above procedure is denoted by $\bm{\mathcal{K}}_{T-1}$.
The outer approximation of the expected future cost can then be defined as
\begin{flalign}
	&\mathfrak {V}_T:=\max_{k\in\bm{\mathcal{K}}^{\rm }_{T-1}}\{\bm{g}_{T,k}^{\top}{\bm{x}}_{T-1}+\beta_{T,k}\}.&
\end{flalign}
Given a trial solution $\overline{\bm{x}}_{T-2}$, a The subproblem in the $t=T-1$ stage can be rewritten as follows:
\begin{subequations}
	\begin{flalign}
		\label{eq:HJB_backwardT-1}
		&\min_{\bm{x}_{T-1}\ge0,\mathfrak{V}_{T}\in\mathbb{R}, \bm{A}^n_{T-1}\bm{x}_{T-1}=-\bm{B}^n_{T-1}\overline{\bm{x}}_{T-2}+\bm{b}^n_{T-1}}(\bm{c}_{T-1}^{n })^{\top}\bm{x}_{T-1} + \mathfrak{V}_{T}&\forall n=1,...,N_{T-1},\\
		\label{eq:BendersCutting_T-1}
		&\mathfrak{V}_{T}\ge(\bm{x}_{T-1})\bm{g}_{T,k}^{\top}+\beta_{T,k}&\forall k\in\bm{\mathcal{K}}_{T-1}.
	\end{flalign}
\end{subequations}
If the optimal solution of Eq.(\ref{eq:HJB_backwardT-1}) is well-defined with respect to almost all $\bm{\xi}_{T}$ realizations, the Benders cut can be computed using the same procedure as in $t=T$ problem.
Applying the same computational procedure to $t=T-2,...,1$, we obtain an outer approximation of the future value function at each stage. 

The backward step finally solves the subproblem at the initial step $t=0$:
\begin{subequations}
	\label{eq:single_backward0}
	\begin{flalign}
		\label{eq:single_HJB_backward0}
		&\underline{\Psi}=\min_{\bm{x}_{1}\ge0, \mathfrak{V}_{2}\in\mathbb{R}, \bm{A}_{1}\bm{x}_{1}=-\bm{B}_{1}\bm{x}_{0}+\bm{b}_{1}}\bm{c}_{1}^{\top}\bm{x}_{1} + \mathfrak{V}_{2},&\\
		\label{eq:single_BendersCutting0}
		&\mathfrak{V}_{2}\ge\bm{x}_{1}\bm{g}_{2,k}^{\top}+\beta_{2,k}&\forall k\in\bm{\mathcal{K}}_1.
	\end{flalign}
\end{subequations}
Eq.(\ref{eq:single_backward0}) provides a lower bound on the optimal value of the SSA problem.
Since each iteration of the backward step adds a new Benders cut to Eq.(\ref{eq:single_BendersCutting0}), we can guarantee that the optimal value of Eq.(\ref{eq:single_backward0}) is non-decreasing.
The pseudo-code for the backward step is shown in Algorithm \ref{alg:single_Backward}.

\begin{figure}[!t]
	\begin{algorithm}[H]
		\caption{Backward pass of the SDDP algorithm for $M$ trial solutions}
		\label{alg:single_Backward}
		\begin{algorithmic}[]
			\setlength{\itemsep}{0pt} \setlength{\parskip}{0pt}
			\REQUIRE $\{\overline{\bm{x}}_t^m\}_{t=0,...,T-1, m=1,...,M}$
			\FOR{$t=T\rightarrow2$}
			\FOR{$m=1\rightarrow M$}
			\FOR{$n=1\rightarrow N_t$}
			\STATE Solve $t$-th stage problem $(\mathfrak{V}_{T+1}=0, \underline{V}_T(\overline{\bm{x}}_{T-1}^m,\bm{\tilde{\xi}}^n_T)={V}_T(\overline{\bm{x}}_{T-1}^m,\bm{\tilde{\xi}}^n_T),~ \underline{\mathscr{V}}_T(\overline{\bm{x}}_{T-1}^m)=\mathscr{V}_T(\overline{\bm{x}}_{T-1}^m)):$
			\STATE $\left(\underline{V}_t(\overline{\bm{x}}_{t-1}^m,\bm{\tilde{\xi}}^n_t), \bm{\pi}_t^{n,m}\right)\leftarrow$
			\STATE $\min_{\bm{x}_t\ge0, \mathfrak{V}_{t+1}\in\mathbb{R}}\left\{\bm{c}_t^{n\top}\bm{x}_t+\mathfrak{V}_{t+1}|\bm{A}^n_t\bm{x}_t=-\bm{B}^n_t\overline{\bm{x}}_{t-1}^m+\bm{b}_t^n,~\mathfrak{V}_{t+1}\ge (\bm{g}_{t+1,k}^m)^{\top}\bm{x}_t+\beta^m_{t+1,k},~\forall k\in\bm{\mathcal{K}}_t^{\rm } \right\}$
			\STATE Compute $\bm{g}_t^{n,m}$: $\bm{g}_t^{n,m}=-(\bm{B}^{n}_t)^{\top}\bm{\pi}_t^{n,m}$
			\ENDFOR
			\STATE Compute $\bm{g}_t^m,~\mathscr{V}_t(\overline{\bm{x}}_{t-1}^m),~{\rm and}~\beta_t^m:$
			\STATE $\bm{g}_t^m=\frac{1}{N_t}\sum_{n=1}^{N_t}\bm{g}_t^{n,m},$ 
			\STATE $\mathscr{V}_t(\overline{\bm{x}}_{t-1}^m)=\frac{1}{N_t}\sum_{n=1}^{N_t}\underline{V}_t(\overline{\bm{x}}_{t-1}^m,\bm{\tilde{\xi}}^n_t),$
			\STATE $\beta_t^m=\mathscr{V}_t(\overline{\bm{x}}_{t-1}^m)-(\bm{g}_t^m)^{\top}\overline{\bm{x}}_{t-1}^{m}.$
			\STATE Add the new cut $\mathfrak{V}_t\ge (\bm{g}_t^m)^{\top}\bm{x}_{t-1}+\beta^m_{t}$ to all $(t-1)$-th stage problems
			\ENDFOR
			\ENDFOR
			\STATE Update the lower bound: Solve 1st stage problem
			\STATE $\underline{\Psi}\leftarrow\min_{\bm{x}_1\ge0, \mathfrak{V}_{2}\in\mathbb{R}}\left\{\bm{c}_1^{\top}\bm{x}_1+\mathfrak{V}_{2}|\bm{A}_1\bm{x}_1=-\bm{B}_{1}\bm{x}_{0}+\bm{b}_1,~\mathfrak{V}_{2}\ge \bm{x}^{\top}_1\bm{g}_{2,k}+\beta_{2,k},~\forall k\in\bm{\mathcal{K}}^{\rm }_1 \right\}$
		\end{algorithmic}
	\end{algorithm}
\end{figure}

\subsubsection{Forward Step}
The forward step provides a trial solution to the backward step and evaluates the upper bound of the SSA problem.
The forward step solves the subproblem forward from the initial step after sampling $M$ realizations of the stochastic process.
Using the Benders cut set computed in the backward step, the subproblem at the time step $t$ is rewritten as
\begin{subequations}
	\label{eq:forward}
	\begin{flalign}
		\label{eq:HJB_Forward}
		&\overline{\bm{x}}_t^m=\argmin_{\bm{x}_t\ge0, \mathfrak{V}_{t+1}\in\mathbb{R}, \bm{A}^m_t\bm{x}_t=-\bm{B}^m_t\overline{\bm{x}}_{t-1}^m+\bm{b}^m_t}(\bm{c}^{m}_t)^{\top}\bm{x}_t + \mathfrak{V}_{t+1}&\forall t,m=1,...,M,\\
		\label{eq:BendersCutting0_Forward}
		&\mathfrak{V}_{t+1}\ge\bm{g}_{t+1,k}^{\top}\bm{x}_{t}+\beta_{t+1,k}{\rm~where~}\mathfrak{V}_{T+1}=0&\forall t,k\in\bm{\mathcal{K}}_t.
	\end{flalign}
\end{subequations}
Eq.(\ref{eq:forward}) can be computed in parallel for each sampled realization.
The pseudo-code for the forward step is shown in Algorithm \ref{alg:Forward}.

\begin{figure}[!t]
	\begin{algorithm}[H]
		\caption{Forward pass of the SDDP algorithm for $M$ sample paths}
		\label{alg:Forward}
		\begin{algorithmic}[]
			\setlength{\itemsep}{0pt} \setlength{\parskip}{0pt}
			\REQUIRE Finite lower bounds on $\{\mathfrak{V}_t\}_{t=2,...,T}$
			\FOR{$m=1\rightarrow M$}
			\STATE Solve 1st stage problem:
			\STATE $\overline{\bm{x}}_1=\argmin_{\bm{x}_1,\mathfrak{V}_2\in\mathbb{R}}\{\bm{c}_1^{\top}\bm{x}_1+\mathfrak{V}_{2}|\bm{A}_1\bm{x}_1=-\bm{B}_{1}\bm{x}_{0}+\bm{b}_1,~\mathfrak{V}_{2}\ge \bm{g}_{2,k}^{\top}\bm{x}_1+\beta_{2,k},~\forall k\in\bm{\mathcal{K}}^{}_1\}$
			\STATE Generate sample path $m$: $\{\bm{\tilde{\xi}}_2^m,...,\bm{\tilde{\xi}}_T^{m}\}$
			\FOR{$t=2\rightarrow T$}
			\STATE Solve $t$-th stage problem ($\mathfrak{V}_{T+1}=0$):
			\STATE $\overline{\bm{x}}_t^m=\argmin_{\bm{x}_t,\mathfrak{V}_{t+1}\in\mathbb{R}}\{\bm{c}_t^{\top}\bm{x}_t+\mathfrak{V}_{t+1}|\bm{A}^m_t\bm{x}_t=-\bm{B}^m_t\overline{\bm{x}}_{t-1}+\bm{b}^m_t,~\mathfrak{V}_{t+1}\ge \bm{g}_{t+1,k}^{\top}\bm{x}_t+\beta_{t+1,k},~\forall k\in\bm{\mathcal{K}}^{\rm }_t\}$
			\ENDFOR
			\STATE Compute the objective value:
			\STATE ${\Psi}^m\leftarrow \bm{c}_1^{\top}\overline{\bm{x}}_1+\sum_{t=2}^T(\bm{c}_t^m)^{\top}\overline{\bm{x}}_t^m$
			\ENDFOR
			\STATE Compute the sample average $\mu$ and SD $\sigma$:
			\STATE $\mu=\frac{1}{M}\sum_{m=1}^M{\Psi}^m$ and $\sigma^2=\frac{1}{M-1}\sum_{m=1}^M({\Psi}^m-\mu)^2$
			\STATE Update the upper bound: $\overline{\Psi}\leftarrow\mu+z_{\alpha}\sigma\sqrt{M}$
		\end{algorithmic}
	\end{algorithm}
\end{figure}

\subsection{Reformulation}

The proposed [SAV-MSSLP] can satisfy most of the assumption discussed in the previous section, which are necessary to obtain an optimal solution with SDDP algorithm.
Unfortunately, [SAV-MSSLP] does not satisfy relatively complete recourses due to demand attraction constraints (\ref{eq:Destination_Demand})(\ref{eq:Destination_Demand_hat}) and supply attraction constraints (\ref{eq:SAV_return})(\ref{eq:SAV_return0}).
A promising solution is to relax the constraints that violate the assumptions and add them to the objective function as penalty functions. 
Then, using a sufficiently large constant $\alpha_P$, objective function in Eq.(\ref{eq:ObjectiveFunction}) is rewritten as follows:
\begin{flalign}
	\label{eq:NewObjectiveFunction}
	&F^t=\begin{cases}
		\alpha_C C & {\rm if}~~t = 0\\
		\alpha_N N & {\rm if}~~t = 1\\
		\alpha_T T^t + \alpha_D D^t + \alpha_P\left(\sum_{k,s}P_{s}^{k,t} + \sum_iP_i^t\right)& {\rm otherwise} 
	\end{cases}&\forall t\in\mathcal{T}_0,\\  
	&P_{s}^{k,t} = \begin{cases}
		\left(\sum_{r\in\mathcal{R}}\xi_{rs}^k - Q_s^{k,t}\right) + \left(\sum_{r\in\mathcal{R}}\hat{\xi}_{rs}^k - \hat{Q}_s^{k,t}\right)& {\rm if~}t\in[T^k,...,T]\\
		0& {\rm otherwise~}
	\end{cases}&\forall k,s,t\in\mathcal{T},\\
	&P_i^t = \begin{cases}
		x_{i0}^t + \hat{x}_{i0}^t & {\rm if~}t\neq T\\
		0& {\rm otherwise~}
	\end{cases}&\forall i,t\in\mathcal{T},\\
	&x_{i0}^t\ge0&\forall i,t\in\mathcal{T},\\
	&\hat{x}_{i0}^t\ge0&\forall i,t\in\mathcal{T}.
\end{flalign}

\section{Numerical Examples}
In this section, numerical examples examine the properties and performance of SAV system design with pre-booked and on-demand requests.
Our numerical experiments are divided into three parts.
We first verified the convergence properties of the SDDP.
Second, we clarified the effect of pre-booking options and dedicated SAVs through sensitivity analysis. 
The sensitivity analysis additionally explored the Pareto-efficiency of the proposed system.
Third, comparing the pre-booked SAV system with on-demand SAV system demonstrated the importance of decision-making considering pre-booking options.

\subsection{Numerical Settings}
The following parameter values were considered during all numerical examples:
\begin{itemize}
    \setlength{\itemsep}{0pt} \setlength{\parskip}{0pt}
    \item We focus one-dimensional urban network consisting of five cities.
    \item The length of planning horizon is given by $T=16$.
    \item The latest arrival time of traveler is given by $T^k=10$.
    \item Infrastructure unit costs are given by $c_i=1~~\forall i$ and $c_{ij}=4~~\forall ij$.
    \item The free flow travel time and waiting time are given by $t_{ij}=1$ and $t_{ii}=1$.
    \item The length of link is also given by $d_{ij}=1$.
    \item The maximum and minimum allowable value of $\mu^t_{ij}$ are $\mu_{ij}^{\rm max}=80$ and $\mu_{ij}^{\rm min}=20$.
    \item The weighted parameters are given in three types: $(\alpha_T,\alpha_D,\alpha_N,\alpha_C)=(10,1,1,1)$, $ (1,1,1,1)$, or $(1,10,1,1)$. The base case is $(\alpha_T,\alpha_D,\alpha_N,\alpha_C)=(10,1,1,1)$.
    \item The vehicle carrying capacity are given in two types: $\rho=3$ or $\rho=4$. The base case is $\rho=3$.
    \item The probability distribution of traffic demand pattern is randomly generated so that the expected sum of traffic demand within the planning horizon is 1,000.
    Then, the expected total traffic demand is allocated into that of pre-booked and on-demand requests based on an exogenous booking rate.
    The booking rate is given in five types: 0.0, 0.25, 0.5, 0.75, and 1.0. The base booking rate is 0.5.
    \item All traffic demand follows a uniform distribution $[0.8U, 1.2U]$, where $U$ is the expected value.
    The departure time step is given by $k=2$ or $k=5$.
    \item The Sample Average Approximation problem is constructed from 1,000 samples.
\end{itemize}

\subsection{Results and Discussion}

\subsubsection{Validation}

\begin{figure}[!ht]
    \centering
    \begin{subfigure}{0.48\textwidth}
        \centering
        \includegraphics[width=\textwidth]{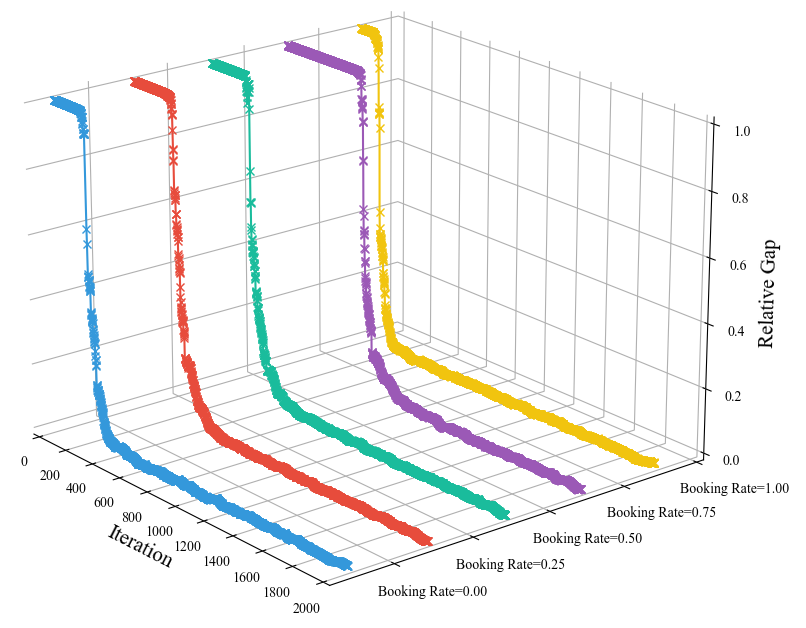}
        \caption{$(\alpha_t = 10, \alpha_t = 1)$}
        \label{fig:convergence1}
    \end{subfigure}
    \hfill
    \begin{subfigure}{0.48\textwidth}
        \centering
        \includegraphics[width=\textwidth]{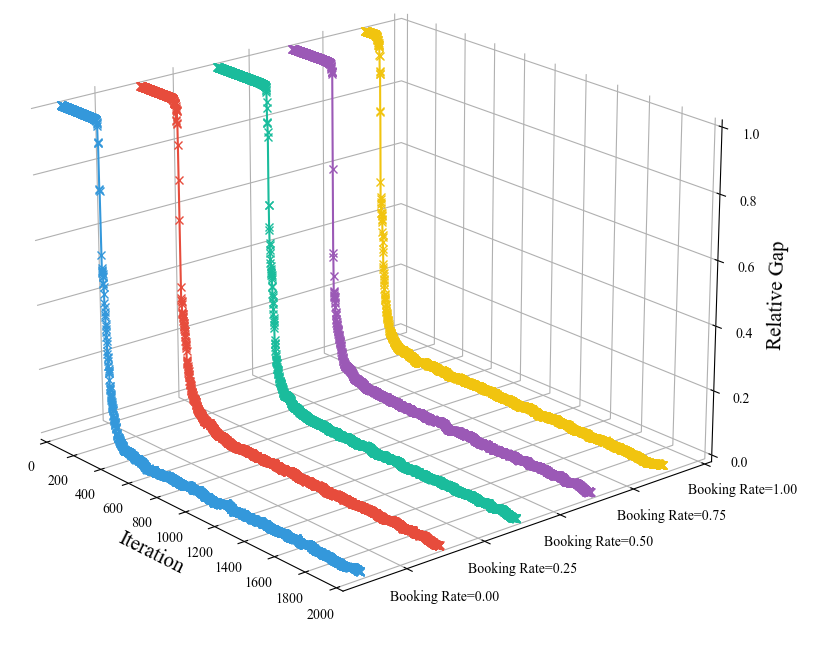}
        \caption{$(\alpha_t = 1, \alpha_t = 10)$}
        \label{fig:convergence2}
    \end{subfigure}
    \caption{Convergence properties}
    \label{fig:convergence}
\end{figure}

Figure~\ref{fig:convergence} shows convergence properties of the solution to the [SAV-MSSLP] by SDDP.
The x-axis depicts the number of iterations, the y-axis depicts parameter settings, 
and the z-axis depicts relative gaps.
The left and right in Figure~\ref{fig:convergence} are convergence properties when weighted parameters are given by $(\alpha_T=10, \alpha_D=1)$ and $(\alpha_T=1, \alpha_D=10)$, respectively.
Figure~\ref{fig:convergence} shows that after 1,000 iterations, the results for all parameter settings converge with a relative gap of almost 0.0, indicating that the obtained solution is close enough to the global optimum.

\subsubsection{Sensitivity Analysis on SAV System Design}

\begin{figure}[!ht]
  \centering
  \includegraphics[width=0.9\textwidth]{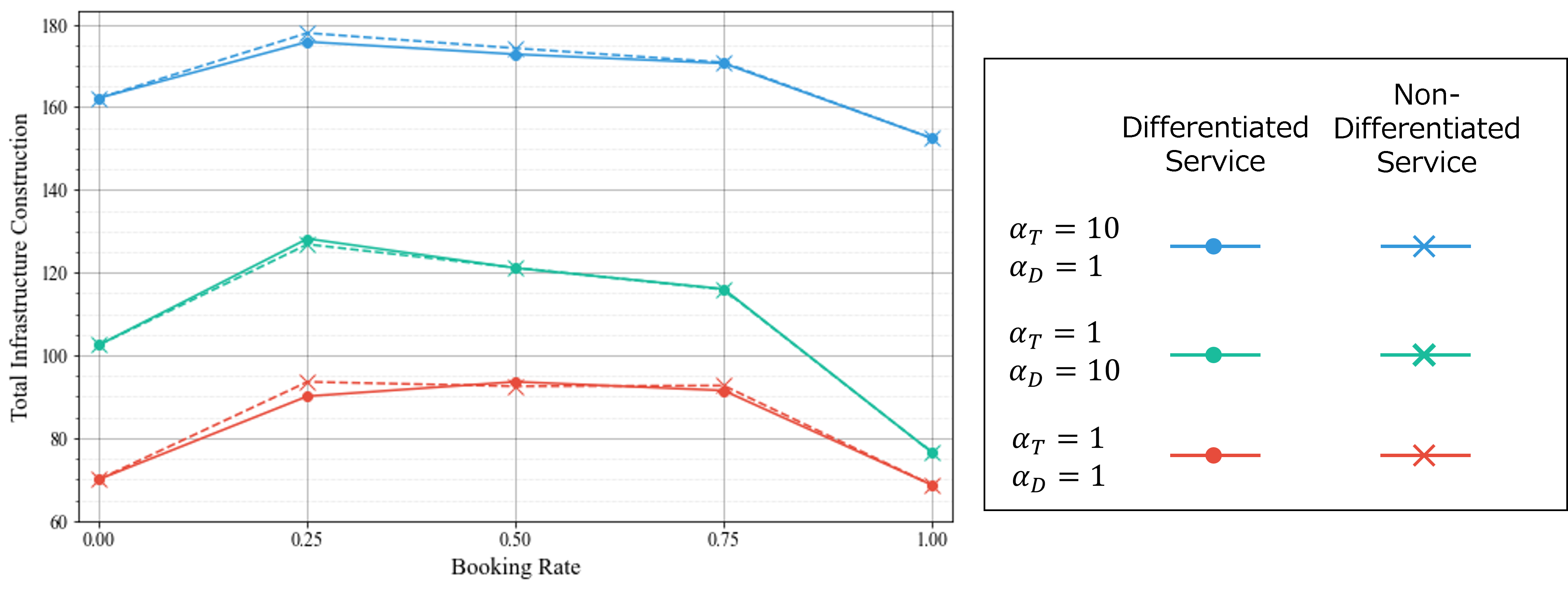}
  \caption{Sensitivity analysis of total infrastructure construction.} \label{fig:construction}
\end{figure}

\begin{figure}[!ht]
  \centering
  \includegraphics[width=0.9\textwidth]{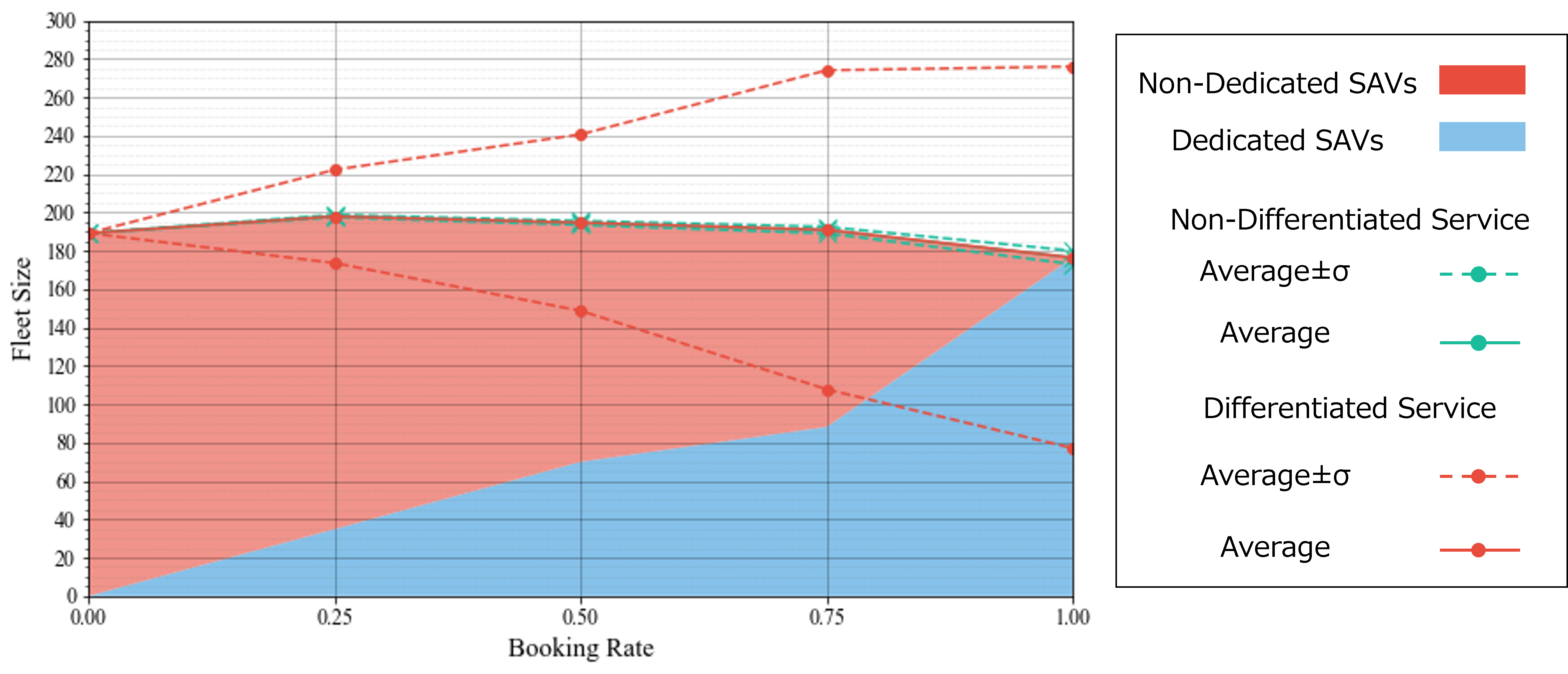}
  \caption{Sensitivity analysis of fleet size.} \label{fig:fleet}
\end{figure}

Figures~\ref{fig:construction} and \ref{fig:fleet} explore the sensitivity of infrastructure planning and fleet sizing strategies with the introduction of pre-booking options and dedicated SAVs.
The horizontal axes in both Figures mark booking rate.
The vertical axes in Figures~\ref{fig:construction} and \ref{fig:fleet} show total infrastructure construction and the number of SAVs, respectively.
The circle and cross marks in Figure~\ref{fig:construction} represent the total infrastructure construction in SAV systems with and without dedicated SAVs, respectively.
Similarly, in Figure \ref{fig:fleet}, the red and green lines represent the fleet size in SAV systems with and without dedicated SAVs, respectively.

From Figures~\ref{fig:construction} and \ref{fig:fleet}, we can see that the introduction of the pre-booking option tends to increase infrastructure construction and the number of vehicles, versus the introduction of dedicated vehicles, which does not affect planning and strategic decisions significantly.
This suggests that faster acquisition of reliable information affects the benefits of investing in vehicles and infrastructure resources.
In contrast, Figures~\ref{fig:construction} and \ref{fig:fleet} illustrate that the higher the booking rate, the lower the investment in vehicles and infrastructure.
This suggests that pre-booking options can mitigate demand uncertainty, resulting in decreased number of vehicles required.
Comparing the fleet size variation of SAV systems with and without dedicated vehicles in Figure~\ref{fig:fleet} shows that the fleet sizing strategy fluctuates significantly with the introduction of dedicated SAVs.
Moreover, as the number of pre-booked travelers increases, the variability in fleet size increases.
It can be concluded that dedicated vehicles are more sensitive to the fluctuations of pre-booked travelers because they do not provide mobility services to on-demand travelers.

\subsubsection{Performance Evaluation}

\begin{figure}[!ht]
  \centering
  \includegraphics[width=0.9\textwidth]{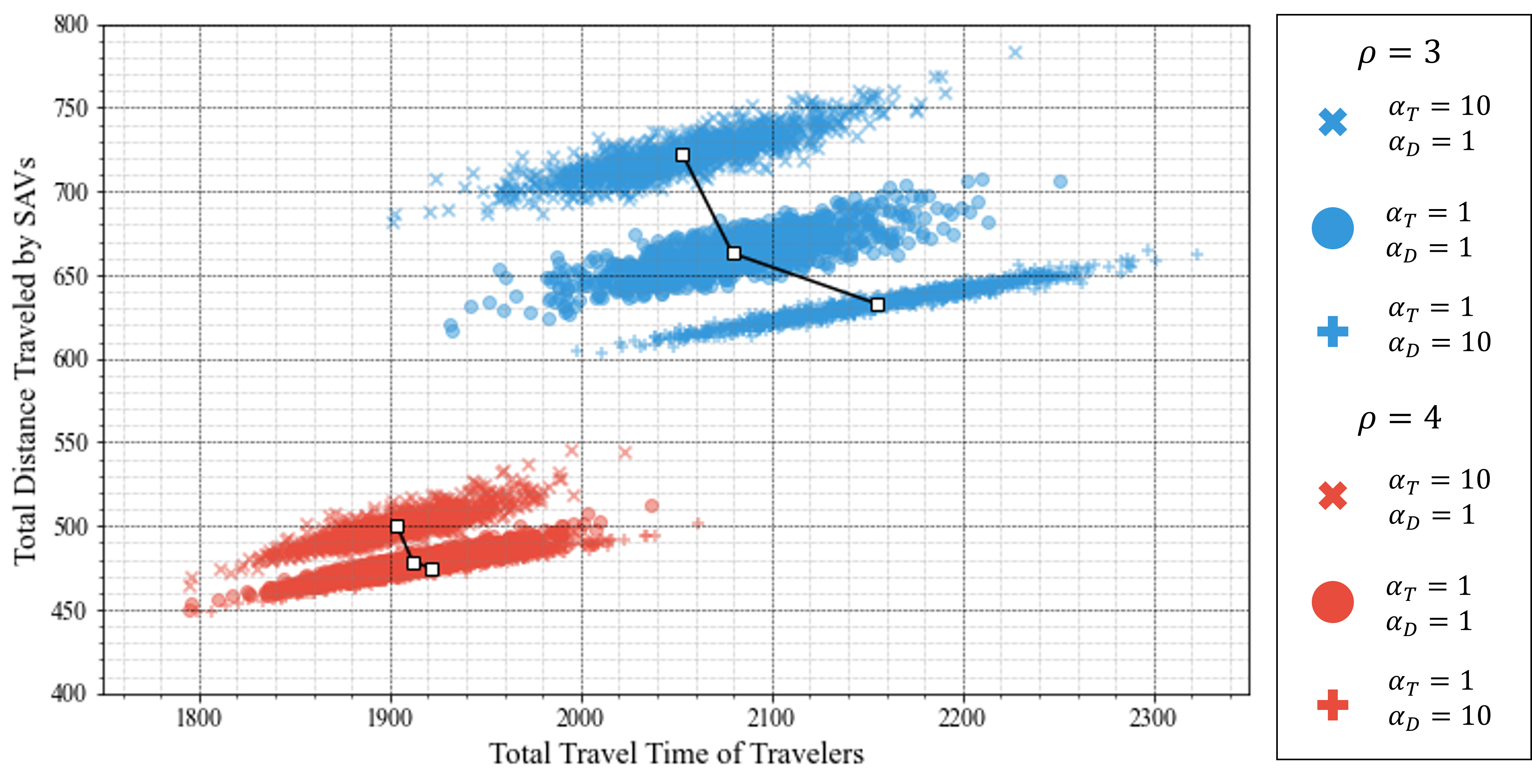}
  \caption{Pareto-efficiency by expanding vehicle carrying capacity.} \label{fig:parate1-2}
\end{figure}

Figure~\ref{fig:parate1-2} shows the approximated Pareto frontier in terms of expected performance.
The actual Pareto frontier is a four-dimensional polyhedron, which is impossible to illustrate.
To illustrate a feature of the Pareto frontier, we draw a two-dimensional relation (i.e., a cross-section of the four-dimensional Pareto frontier) between total travel time, $T$, and total distance traveled by SAVs, $D$.
The horizontal and vertical axes depict the total travel time of travelers and the total distance traveled by SAVs, respectively.
The cross, circle, and plus marks represent performance indexes when weighted parameters are given by $(\alpha_T=10, \alpha_D=1)$, $(\alpha_T=1, \alpha_D=1)$, and $(\alpha_T=1, \alpha_D=10)$, respectively 
The blue and red marks represent performance indexes when vehicle carrying capacities are given by $\rho=3$ and $\rho=4$, respectively.

Figure~\ref{fig:parate1-2} indicates that the weighted sum method can provide a Pareto-efficient solution in terms of expected value.
We also point out that, obviously, increasing vehicle carrying capacity improves the Pareto-efficiency.
This quantitative property always holds in a deterministic framework (\citep{seo2021multi}).

\begin{figure}[!ht]
  \centering
  \includegraphics[width=0.9\textwidth]{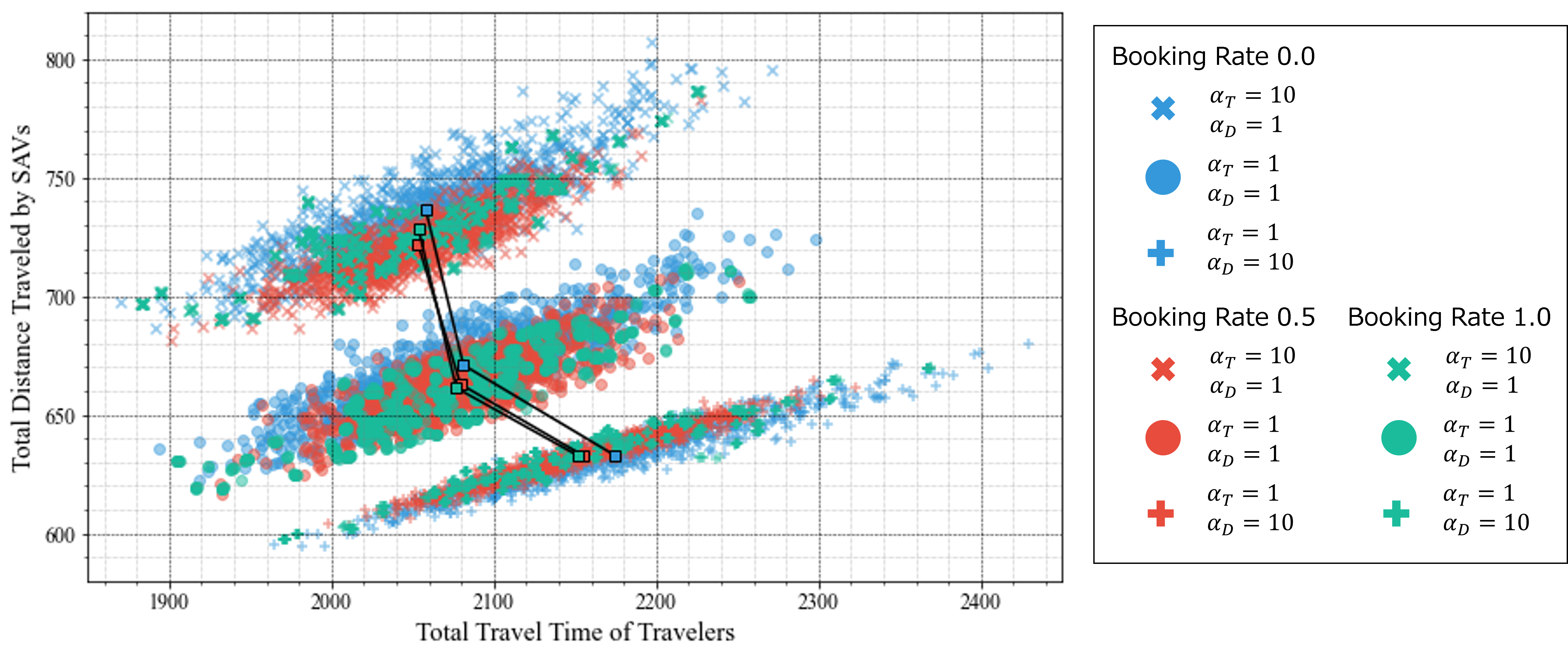}
  \caption{Pareto-efficiency by introducing pre-booking options.} \label{fig:parate2}
\end{figure}

\begin{figure}[!ht]
  \centering
  \includegraphics[width=0.9\textwidth]{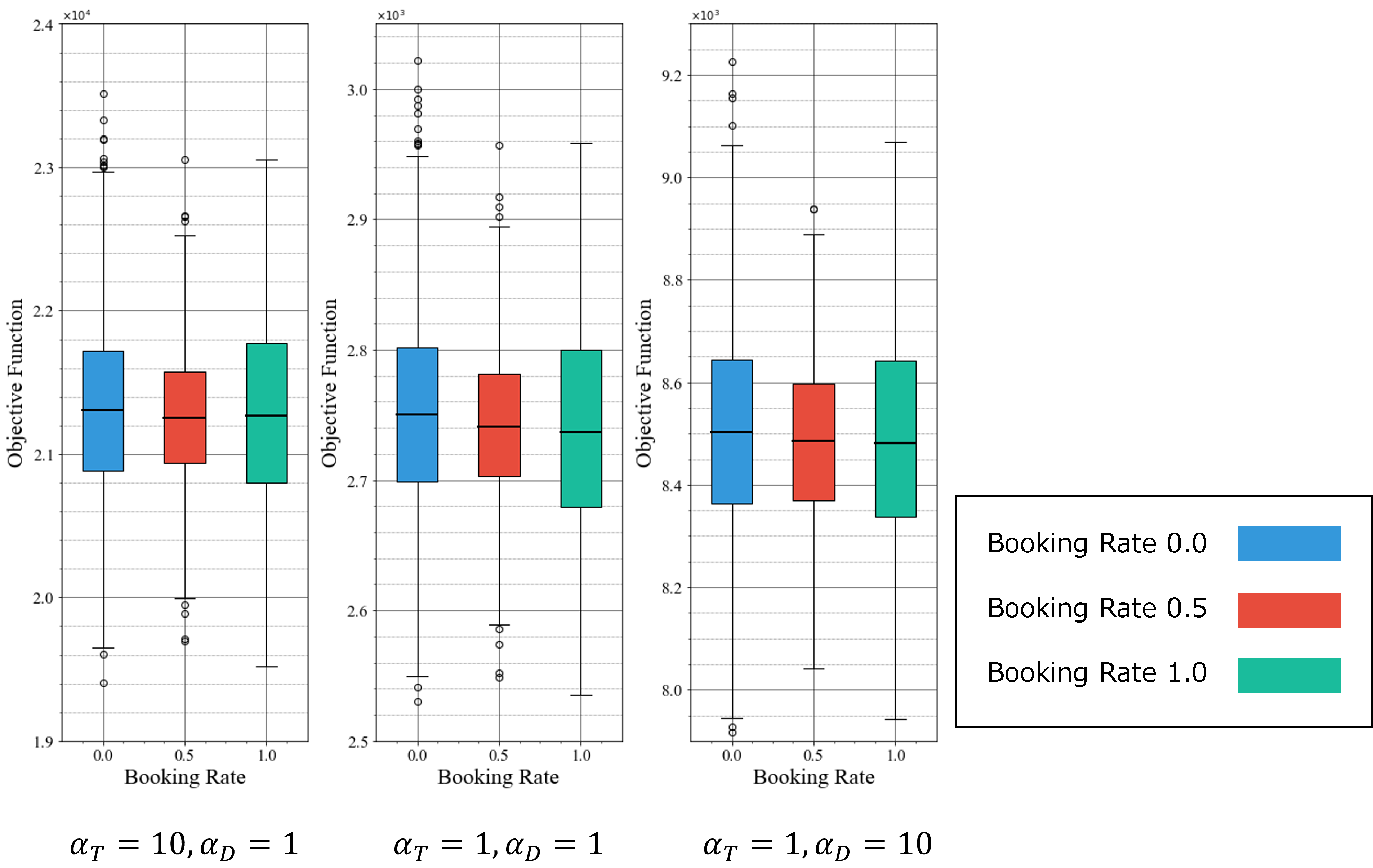}
  \caption{Objective function by introducing pre-booking options.} \label{fig:parate2_BOX}
\end{figure}

Figures~\ref{fig:parate2} and \ref{fig:parate2_BOX} show the impact on Pareto-efficiency and system performance by introducing pre-booking options, respectively.
The representation of the vertical and horizontal axes and marks in Figure~\ref{fig:parate2} is the same as in Figure~\ref{fig:parate1-2}.
The blue, red, and green marks represent performance indexes when booking rates are given by $0.0$, $0.5$, and $1.0$, respectively.
Figure \ref{fig:parate2_BOX} are box plots of the objective function.
The left in Figure \ref{fig:parate2_BOX} shows the objective function when the weighted parameters are given as $(\alpha_T=10, \alpha_D=1)$, $(\alpha_T=1, \alpha_D=1)$, and $(\alpha_T=1, \alpha_D=10)$, in that order.
The colors in the box plots correspond to the marks in Figure \ref{fig:parate2}.

Figure~\ref{fig:parate2} indicates the Pareto-improvement by introducing pre-booking options in terms of expected performance.
We note that this Pareto improvement is not significant.
Figures~\ref{fig:parate2} and \ref{fig:parate2_BOX} also suggest that the increase in the booking rate does not contribute much to improving system performances.

\begin{figure}[!ht]
  \centering
  \includegraphics[width=0.9\textwidth]{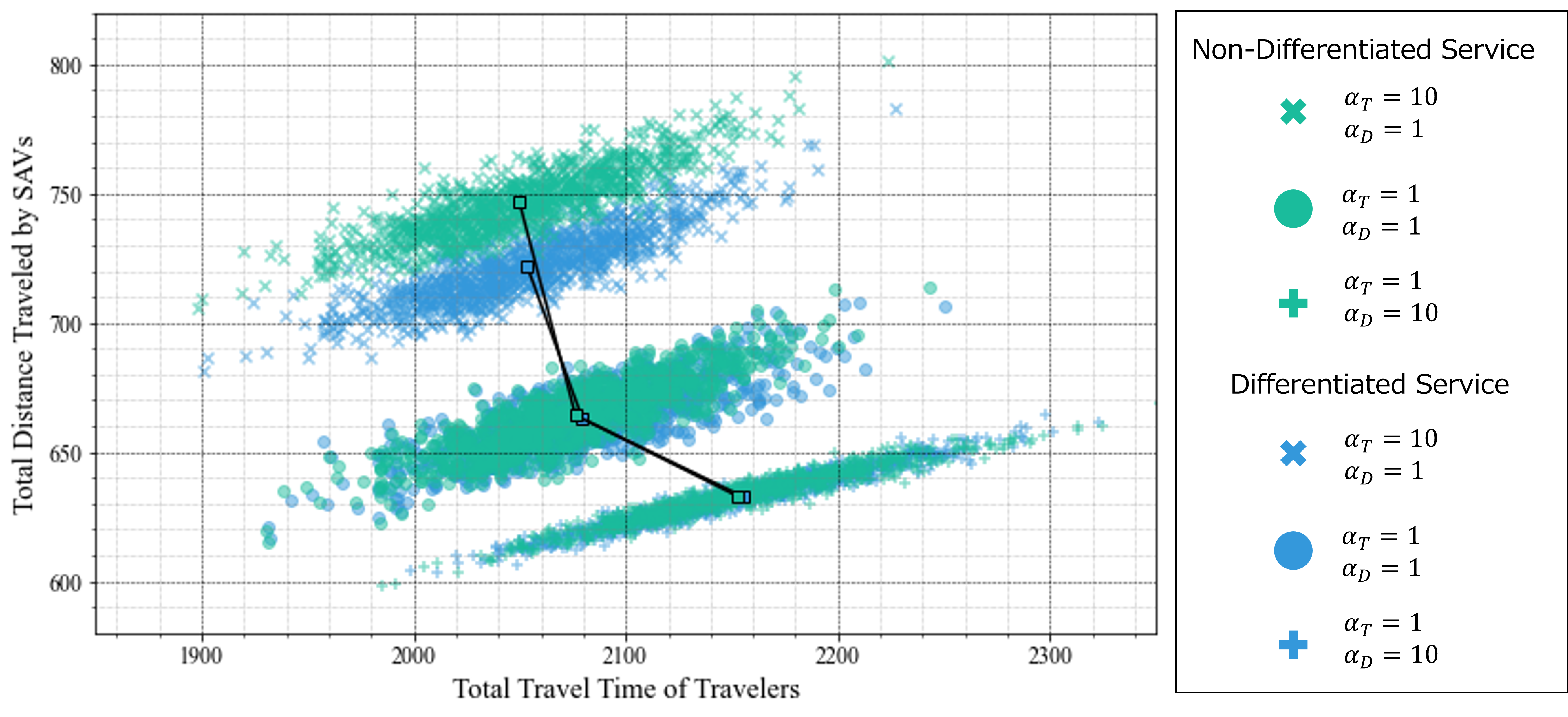}
  \caption{Pareto-efficiency by introducing dedicated SAVs.} \label{fig:parate3}
\end{figure}

Figure~\ref{fig:parate3} shows the impact on Pareto-efficiency by introducing dedicated SAVs.
The representation of the vertical and horizontal axes and marks in Figure~\ref{fig:parate3} is the same as in Figure~\ref{fig:parate1-2}.
The blue and green marks in Figures~\ref{fig:parate3} represent performance indexes in SAV system with and without dedicated vehicles, respectively.
From Figure~\ref{fig:parate3}, we can see that the system performance does not change significantly by introducing dedicated vehicles.

\begin{figure}[!ht]
  \centering
  \includegraphics[width=0.9\textwidth]{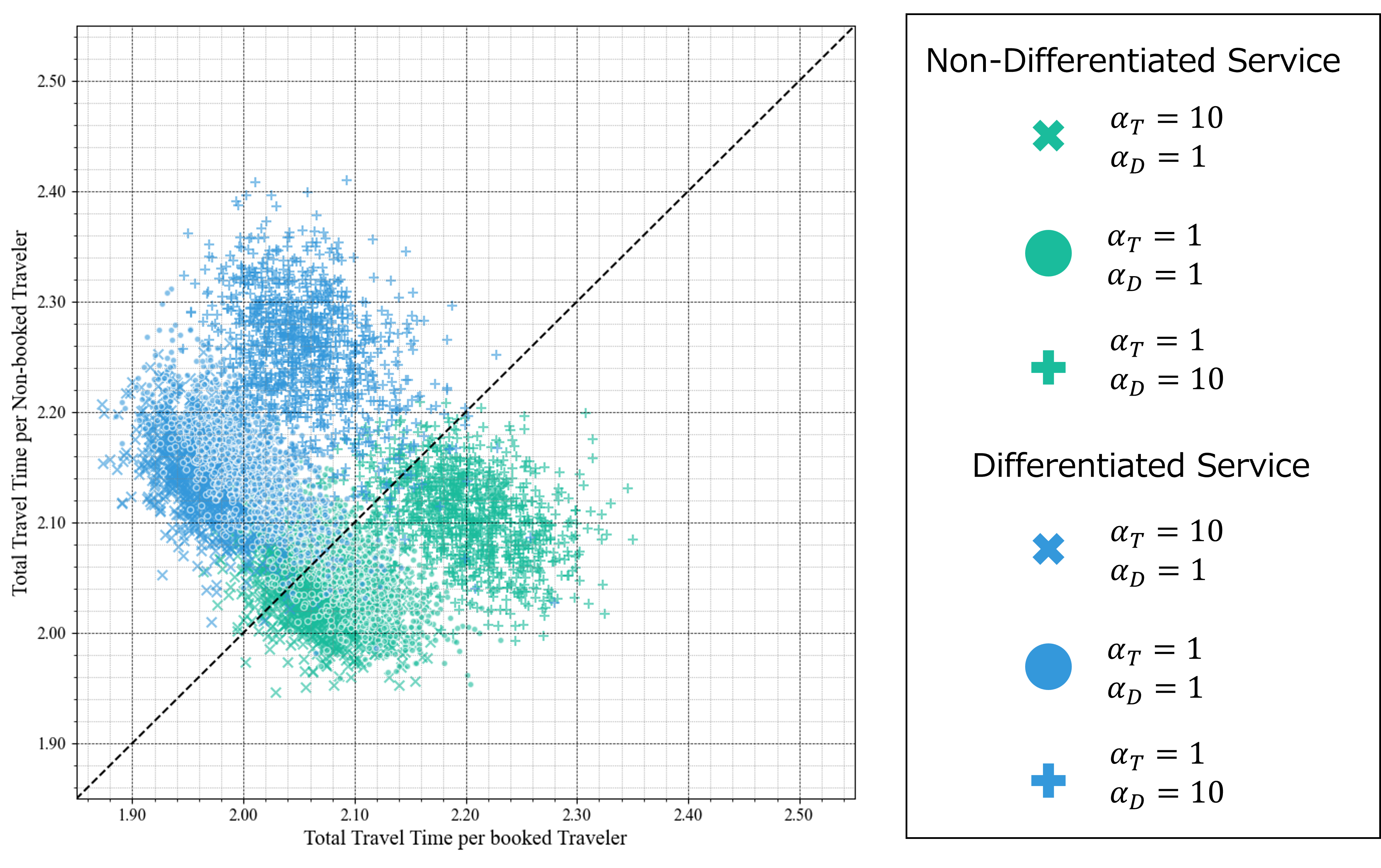}
  \caption{Impact on travel time per traveler by introducing dedicated SAVs.} \label{fig:MarginalTT}
\end{figure}

Figure~\ref{fig:MarginalTT} represents the total travel time per traveler for pre-booked and on-demand travelers.
The horizontal and vertical axes mark that for pre-booked and on-demand travelers, respectively.
The representation of the marks in Figure~\ref{fig:MarginalTT} is the same as in Figure~\ref{fig:parate3}.

Figure~\ref{fig:MarginalTT} points out that dedicated SAVs bring a travel time reduction effect for pre-booked travelers relative to on-demand travelers.
This result suggests that introducing dedicated vehicles can provide an incentive for pre-booking without significantly degrading system performance, thus ensuring the sustainability of the pre-booking system.

\begin{figure}[!ht]
  \centering
  \includegraphics[width=0.9\textwidth]{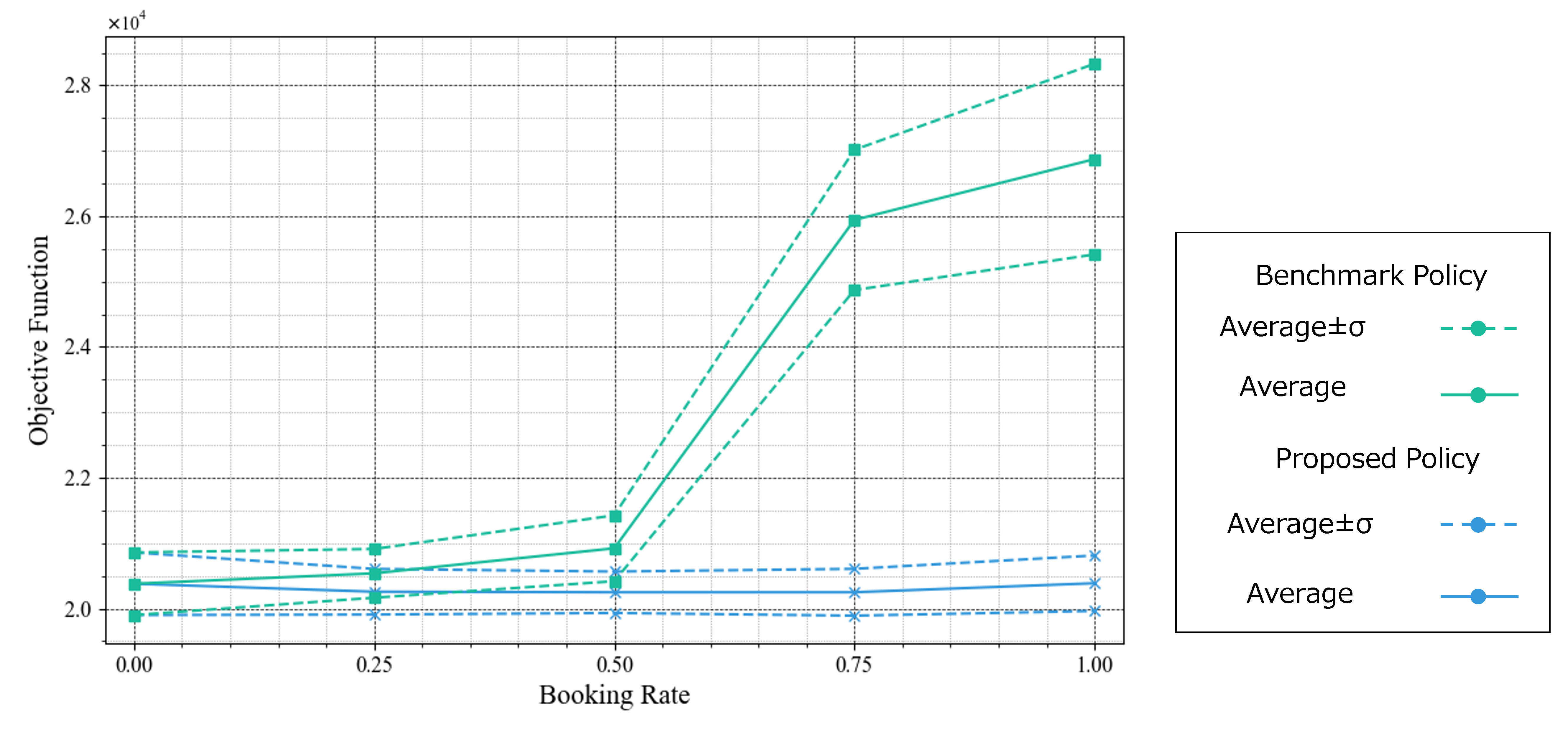}
  \caption{Comparison of performance between systems that assume pre-booking options and do not.} \label{fig:performance}
\end{figure}

Finally, Figure~\ref{fig:performance} evaluates the performance of the SAV system design capturing pre-booking options.
We consider the SAV system discarding pre-booking options as a benchmark policy.
The horizontal and vertical axes depict the booking rate and objective function, respectively.
The blue and green mark and line mean system performance of the proposed system design and benchmark policy.

Figure~\ref{fig:performance} indicates that the proposed system always outperforms the benchmark policy. 
The difference in performance increases as the booking rate increases, and especially when the booking rate exceeds 50\%, the performance gap dramatically increases.
This performance difference would become more noticeable as the scale of the problem increases, due to economies of scale.

\section{Conclusion}

This study presents a joint optimization framework for long-term infrastructure planning, mid-term fleet sizing, and short-term routing and ridesharing matching operations in SAV systems with pre-booked and on-demand requests.
Dynamic traffic assignment approaches enable us to formulate the SAV system design and operation problems as multi-stage stochastic linear programming.
Leveraging the linearity of the problem can help tackle the computational difficulties caused by multi-objective and dynamic stochasticity.

Our numerical examples demonstrate the importance of designing SAV systems with pre-booking options.
Our instances exhibit that, even though the test network is small-scale, the pre-booking option affects the long-term network design and mid-term fleet sizing strategy, resulting in Pareto improvements.
This result is assumed to be more noticeable when faced with city-scale problems.
We also show quantitatively that the introduction of dedicated vehicles to pick up and drop off only pre-booked travelers can incentivize pre-booking with little sacrifice in system performance.
This strategy can ensure the sustainability of the pre-booking system.

It should be noted that this study is designed as an initial step toward developing the SAV system design with pre-booked requests.
First, we need to improve the computational performance so as to make it applicable to urban-scale problems as well.
Exploiting the linearity of the problem, we can extend the framework to solve the problem in parallel per subspace through the Dantzig-Wolfe decomposition.
Furthermore, the complexity of dimensionality in variables would require us to employ an iterative algorithm based on neural networks \citep{dai2021neural, kim2024transformer}.
It is also essential to construct a system that can work robustly even when the demand deviates significantly from the probability distribution.
Especially in long-term network design, it is rare that the precise probability distribution of demand is available.
Stochastic dual-dynamic programming has been extended to a framework that can describe flexible risk measures, such as robust programming \citep{georghiou2019robust} and distributionally robust programming \citep{duque2020distributionally}.
Finally, instead of directly controlling SAVs, a framework for pricing is a future challenge.
In transportation systems that are not shared economy, a pricing system is constructed that is desirable for both users and operators via reservations.
By exploiting the linearity of the problem proposed in this study, we can develop a similar pricing system based on \citet{akamatsu2017tradable}.

\section*{Acknowledgements}

This study was partially supported by JSPS Grant-in-aid (KAKENHI) \#24K01002.

\bibliographystyle{abbrvnat0}
\bibliography{references}

\end{document}